\documentstyle[12pt,fleqn,dina4,amssymb]{article}

\def\x{\times}
\def\ds{\displaystyle}

\def\bigspitem{\par\hangone\textputone}
\def\hangone{\hangindent 3\parindent}
\def\textputone#1{\noindent
	 \hbox to 3\parindent{\hss#1\enspace}\ignorespaces}

\begin{document}
\sloppy

\newcommand{\fl}{\mbox{ {\rm fl}} }
\newcommand{\afl}{\mbox{ {\rm afl}} }
\newcommand{\acl}{\mbox{ {\rm acl}} }
\newcommand{\hcl}{\mbox{ {\rm hcl}} }
\newcommand{\gsh}{\mbox{ {\rm gsh}} }
\newcommand{\sh}{\mbox{ {\rm sh}} }
\newcommand{\GG}{\mathop{\rm I\!\!\! G}\nolimits}
\newcommand{\Id}{\mbox{ {\rm Id}} }
\newcommand{\id}{\mbox{ {\rm id}} }
\newcommand{\gf}{\mbox{ {\rm gf}} }
\newcommand{\cl}{\mbox{ {\rm cl}} }
\newcommand{\dpo}{\mathop{ \to\!\!\!\!\!\to }\nolimits}
\newcommand{\supp}{\mbox{ {\rm supp}} }
\newcommand{\ver}{\mbox{ {\rm ver}} }
\newcommand{\edg}{\mbox{ {\rm edg}} }
\newcommand{\diam}{\mbox{ {\rm diam}} }
\newcommand{\FP}{\mbox{ {\rm FP}} }
\newcommand{\FH}{\mbox{ {\rm FH}} }
\newcommand{\Res}{\mbox{ {\rm Res}} }
\newcommand{\dpm}{\mathop{\ \succ\!\!\to\!\!\!\!\!\to }\nolimits}
\newcommand{\pfm}{\mathop{\ \succ\!\!\to}\nolimits}
\newcommand{\NN}{\mathop{\rm I\! N}\nolimits}
\newcommand{\QQ}{\mathop{\rm I\!\!\! Q}\nolimits}
\newcommand{\RR}{\mathop{\rm I\! R}\nolimits}
\newcommand{\ZZ}{\mathop{\sl Z\!\! Z}\nolimits}
\newcommand{\BB}{\mathop{\rm I\! B}\nolimits}
\newcommand{\CC}{\mathop{\rm C\!\!\! I}\nolimits}
\newcommand{\FF}{\mathop{\rm I\! F}\nolimits}
\newcommand{\Hom}{\mathop{\rm Hom}\nolimits}
\newcommand{\HNN}{\mathop{\rm Hnn}\nolimits}
\newcommand{\Isom}{\mathop{\rm Isom}\nolimits}
\newcommand{\EE}{\mathop{\rm I\! E}\nolimits}
\newcommand{\HH}{\mathop{\rm I\! H}\nolimits}
\newcommand{\SL}{\mathop{ {\rm SL}} }
\newcommand{\GL}{\mathop{ {\rm GL}} }
\renewcommand{\thefootnote}{\fnsymbol{footnote}}

\centerline{\Large\bf Connectivity properties of group}
\centerline{\Large\bf actions on non-positively curved spaces I:}
\centerline{\Large\bf Controlled connectivity and openness results}
~\\

\centerline{\large by}
~\\

\centerline{\large Robert Bieri\footnote{The first-named author was
supported in part by a grant from the Deutsche Forschungsgemeinschaft.}
and Ross Geoghegan\footnote{The second-named
author was supported in part by a grant from the National Science Foundation.}}
~\\

~\\

\section{Introduction}
\renewcommand{\thefootnote}{\arabic{footnote}}
\setcounter{footnote}{0}

A thorough outline of this paper is given in \S2.  In this introduction we
give a quick indication of what the paper and (to some extent)
its sequel [BG$_{\mbox{II}}$] are about.

Let $G$ be a group 
of type\footnote{$G$ has type $F_n$ if there is a
  $K(G,1)$-complex with finite $n$-skeleton. All groups have type $F_0,
  F_1$ is ``finitely generated'', $F_2$ is ``finitely presented'', etc.}
$F_n$, let $(M,d)$ be a simply connected proper ``non-positively curved''
(i.e. CAT(0)) metric space, and let $\rho : G \to \Isom(M)$ be an action of
$G$ on $M$ by isometries.  
We introduce a topological
property called (uniform) ``controlled $(n-1)$-connectedness'', abbreviated
to $CC^{n-1}$, which such an action $\rho$ may or may not possess.
This property is extracted and generalized from work on the ``geometric
invariants $\Sigma^n(G)$ of the group $G$ (e.g. [BS 80], [BNS 87], [BRe 88],
[Re 88, 89], [BS]) which applies in the special case where
$M = \EE^m$ and $\rho$ is an action by translations.
Experience with this special case suggests that $CC^{n-1}$ is a
fundamental property, and we are encouraged in this belief by two 
instances in which $CC^{n-1}$ turns out to be equivalent to something familiar:
\begin{itemize}
\item The isometric action $\rho$ is uniformly $CC^{-1}$ if and only if 
$\rho$ is cocompact. See Proposition 3.2.
\item When the isometric action $\rho$ has discrete orbits, 
$\rho$
is uniformly $CC^{n-1}$ if and only if $\rho$ is cocompact and the stabilizer
  $G_a$ of some (equivalently, any) point $a \in M$ has type $F_n$. See
  Theorem A.
\end{itemize}
Moreover, we have an openness theorem:  
\begin{itemize}
\item The property uniformly $CC^{n-1}$ is an open condition on the space
of actions $\Hom(G,\Isom(M))$; in other words, if $\rho$ has it so do all
  isometric actions $\rho'$ near $\rho$ in the compact-open topology.
See Theorem B.
\end{itemize}
Of course, as a consequence we obtain new openness results for cocompact actions
in general, and for
cocompact actions with discrete orbits and stabilizers of type $F_n$. The
former is related to results of A.~Weil (see \S 2.5) but is new
in our generality.  Without the added hypothesis of discrete orbits it is
not true that ``cocompact and point stabilizers having type $F_n$'' is an
open condition.  So ``uniformly $CC^{n-1}$'' appears to be the right
concept.

The art of deciding whether a given isometric action $\rho$ of $G$ on
$M$ is $CC^{n-1}$ is still in its infancy.  One approach is to define an
analogous concept ``$\rho$ is $CC^{n-1}$ over $e$'' for the points $e$ ``at
infinity''.  (The set of such points is often called the ``boundary'' of
$M$ and is denoted in these papers by $\partial M$.  With an appropriate
topology, $\partial M$ compactifies $M$ nicely; see \S10.1.)
We define ``$CC^{n-1}$ over $e$'' in the sequel paper [BG$_{\mbox{II}}$]
where we prove the following ``boundary criterion'':
\begin{itemize}
\item{}
If $M$ is almost geodesically complete then $\rho$ is uniformly $CC^{n-1}$
if and only if $\rho$ is $CC^{n-1}$ over each point $e \in \partial M$. See
Theorem H of [BG$_{\mbox{II}}$].  
\end{itemize}
We have used this boundary criterion to work out the $CC^{n-1}$ properties
of the natural action of the group $\SL_2(\ZZ[\frac{1}{m}])$
on the hyperbolic plane
for the case when $m$ is a prime and we have a good idea how to handle the case
when $m$ is an arbitrary natural number, see \S2.8.  We also understand the
situation where $M$ is a locally finite simplicial tree. 

It turns out that whether a given isometric action $\rho$ is $CC^{n-1}$
over $e$ depends in a delicate way upon the point $e\in \partial M$.  Therefore
the subset of $\partial M$,
\[
\Sigma^n(\rho) := \{e\in \partial M
\mid \rho \ \mbox{ is }\ CC^{n-1}\ \mbox{ over }\ e\},
\]
becomes an interesting invariant of the action $\rho$ even (in fact,
particularly) when 
$\Sigma^n(\rho)$ is not all of $\partial M$.  The
study of $\Sigma^n(\rho)$ is pursued in [BG$_{\mbox{II}}$].  

We have said that $\Sigma^n(\rho)$ constitutes 
a far reaching generalization of the ``Geometric Invariants'' $\Sigma^n(G)$
of the group $G$. 
The reader familiar with that literature, in particular 
with [BGr 84], [Me 94, 95, 96, 97], [Geh], [Ko 96], [Be-Br 97], [Bu]
and [MMV 98] which
compute $\Sigma^n(G)$ for metabelian groups, Houghton groups, Borel
subgroups of Chevalley groups over function fields with finite base fields,
certain soluble groups and right angled Coxeter
groups, as well as direct products and graph products,
will know that $\Sigma^n(G)$ is difficult to compute but is a rich
invariant for those groups $G$ which have infinite Abelianization.
Indeed, $\Sigma^n(G)$ is essentially $\Sigma^n(\rho)$ for
the canonical action $\rho$ of $G$ on $G_{ab} \otimes \RR$.  The study of
$\Sigma^n(\rho)$ for arbitrary isometric 
actions $\rho$ on arbitrary CAT(0) spaces $M$
is seen then as the natural non-commutative extension of what has proved to
be a fruitful commutative case.  We go into more detail on this in
[BG$_{\mbox{II}}$].  
~\\

{\bf Acknowledgments:}
G.~Meigniez has independently obtained some of our main results in the special
case $n=1$ and also has some insight for $n > 1$. We have profited from
discussions with him; in particular he pointed out that our proof of
Theorem A$'$ actually proves the stronger Theorem A. Tom Farrell and
Kai-Uwe Bux also gave us useful insights.

\newpage

\section{Outline, Main Results and Examples}

2.1 {\bf Non-positively curved spaces.} We interpret
``non-positively curved'' to mean that $(M,d)$ is a proper CAT(0)
space. In detail: (i) it is a geodesic metric space: this means that an
 isometric copy
of the closed interval $[0,d(a,b)]$ called a {\em geodesic segment}
joins any two points $a,b \in M$; (ii) for any geodesic triangle
$\Delta$ in $M$ with vertices $a,b,c$ let $\Delta'$ denote a triangle in
the Euclidean plane with vertices $a',b',c'$ and corresponding side
lengths of $\Delta'$ and $\Delta$ equal; let $\omega$ and $\omega'$ be geodesic
segments from $b$ to $c$ and from $b'$ to $c'$ respectively; then for
any $0 \leq t \leq d(b,c),\; d(a,\omega(t)) \leq ||a'-\omega'(t)||$; and
(iii) $d$ 
is {\em proper}, i.e. the closed ball $B_r(a)$ around any $a \in M$ of
any radius $r$ is compact.

In a CAT(0) space the geodesic segment from $a$ to $b$ is unique and
varies continuously with $a$ and $b$. This implies that CAT(0) spaces
are contractible.

For one of our results we only require a weaker property:  that the proper
metric space have unique geodesic segments -- in which case we say $M$ is
a {\em 
unique-geodesic metric space}.  Indeed this is enough to imply 
that geodesic segments vary continuously with their end points, hence
contractibility ([BrHa; I(3)]).  

{\bf Examples of CAT(0)-spaces} are Euclidean space $\EE^m$,
hyperbolic space $\HH^m$, locally finite affine buildings, complete
simply connected open Riemannian manifolds of non-positive sectional
curvature, and any finite cartesian product $(\mathop\Pi\limits_i M_i,d)$ of
CAT(0) spaces $(M_i,d_i)$ with $d(a,b) := (\mathop\sum\limits_i
d_i(a_i,b_i)^2)^{\frac{1}{2}}$. 
~\\

2.2 {\bf Controlled connectivity: the definition of $\bf CC^{n-1}$.}
Controlled topology starts with a {\em control function} $h : X \to M$.  In
our case, $X$ will always be a CW complex and $M$ a metric space; we will
add more structure as we go along.  For $a\in M$ and $r > 0$ we denote by
$X_{(a,r)}$ the largest subcomplex of $X$ lying in $h^{-1}(B_r(a))$.  We say
$X$ is {\em controlled} $(n-1)$-{\em connected} $(CC^{n-1})$ over $a$ (with
respect to $h$) if for all $r\geq 0$ and $-1\leq p\leq n-1$ there exists
$\lambda \geq 0$ such that every map $f : S^p \to X_{(a,r)}$ 
extends\footnote{Previous publications on the geometric invariants used the 
terminology ``$X_{(a,r)}$ is essentially $(n-1)$-connected'' rather than 
``$X$ is $CC^{n-1}$''.  This concept also appears in [FePe 95].} to a map
$\tilde f : B^{p+1} \to X_{(a,r+\lambda)}$.  If $X$ is $CC^{n-1}$ over some
$a\in M$ it is easy to see that $X$ is in fact $CC^{n-1}$ over each point of
$M$, so we can speak of $X$ being $CC^{n-1}$ without reference to a point
$a\in M$.  The number $\lambda$ in the definition of ``$CC^{n-1}$ over $a$''
depends on $a$ and on $r$; sometimes we will write $\lambda(a,r)$ to
emphasize this.  We call $\lambda$ a {\em lag}.  

In this paper we are given a non-negative integer $n$ and a group $G$ of
type $F_n$.  We pick an $n$-dimensional 
$(n-1)$-connected\footnote{$(-1)$-connected means ``non-empty'': the sphere 
$S^{-1}$ is empty and has a unique empty map $S^{-1} \to X$. This map extends 
to the ball $B^0$ if and only if the space is non-empty.  Thus 
$(n-1)$-connected always implies non-empty. This may seem pedantic but will 
be useful in [BG$_{\mbox{II}}$].} CW complex $X^n$ on which $G$ acts freely
on the left as a group of cell permuting homeomorphisms with $G\backslash
X^n$ a finite complex.  In other words $X^n$ is the $n$-skeleton of a
contractible free $G$-CW complex $X$ which is chosen so that $X^n$ is
cocompact.  The metric space $M$ will always be a unique-geodesic metric
space---eventually we will be forced to require that it be CAT(0).  Given an
action $\rho : G \to \mbox{Isom}(M)$ of $G$ on $M$ by isometries, we choose
a $G$-equivariant continuous control function $h : X \to M$; this is
possible because $G$ acts freely on $X$ and $M$ is contractible.  

In this context, suppose $X$ is $CC^{n-1}$ over $a$.  The lag $\lambda$ can
be chosen independent of $a$ if and only if the $G$-action on $M$ is
cocompact\footnote{An action of $G$ on $M$ is {\em cocompact} if there is a 
compact subset $K \subset M$ with $GK = M$.}; see \S3.1.   
In \S 3 we will prove that since $G\backslash X$ has finite
$n$-skeleton, the property of $X$ being $CC^{n-1}$ is
independent of the choice of $X$ and of $h$, i.e., is a property of the
action $\rho$. So, if $X$ is $CC^{n-1}$ we will say that $\rho$ is
$CC^{n-1}$. 
~\\

2.3 {\bf The case of discrete orbits.}
Before we state our main results we interpret $CC^{n-1}$ when the $G$-action on
$M$ has discrete orbits.\footnote{An action of $G$ on $M$ {\it has discrete
orbits} if every orbit is a closed discrete subset of $M$.}
~\\

{\bf Theorem A.}  {\em Let $(M,d)$ be a proper unique-geodesic 
metric space, let $G$
be a group of type $F_n$, and let $\rho : G \to {\rm Isom}(M)$ be a 
cocompact action
which has discrete orbits.   Then $\rho$ is
$CC^{n-1}$ if and only if the stabilizer $G_a$ has type $F_n$.}
~\\

It follows that all
stabilizers $G_a, a \in M$, have the same $F_m$-type -- but that does not
come as a surprise since the point stabilizers of a discrete action by
isometries are easily seen to be pairwise commensurable, so one has type
$F_m$ if and only if the other has type $F_m$.

Theorem A is a consequence of the following homotopy version of K.S. Brown's
finiteness criterion [Br 87, Theorem 2.2]:  
~\\

{\bf $F_n$-Criterion.}  {\em Let $H$ be a group, $Y$ a contractible free $H$-CW
complex and $(K_r)_{r\in \RR}$ an increasing filtration of $Y$ by
$H$-subcomplexes so that $Y = \ds{\bigcup_r} K_r$ and each $K_r$ has
cocompact $n$-skeleton.  Then $H$ is of type $F_n$ if and only if $Y$ is
$CC^{n-1}$ with respect\footnote{In \S2.2 we defined $CC^{n-1}$ using
filtrations which came from control functions, but the definition makes
sense with respect to any filtration.} to the filtration $(K_r)$.}
~\\

We discuss the proof of this criterion in \S8.  

Theorem A follows by setting $Y = X$, $K_r = X_{(a,r)}$ and $H = G_a$.
Clearly $X_{(a,r)}$ is a $G_a$-subcomplex.  The remaining part of the proof,
that each $X^n_{(a,r)}$ is cocompact as a $G_a$-complex, is not hard and is
given in \S8.  

A special case of Theorem A is worth noting. If $N = \ker \rho$ we have
a short exact sequence $N \pfm G \dpo Q$ with $Q \leq \Isom(M)$, and
short exact sequences for the stabilizers $N \pfm G_a \dpo Q_a$. If we
replace the assumption that $\rho$ have discrete orbits 
by the stronger assumption
that the induced action of $Q$ on $M$ be properly
discontinuous\footnote{An action of $Q$ on $M$ is {\em properly
    discontinuous} if every point $a \in M$ has a neighbourhood $U$ such
  that $\{q \in Q|qU \cap U \not= \emptyset\}$ is finite (equivalently: if
  the action has discrete orbits and has finite point stabilizers, see Lemma
  8.5).}
 then Theorem A applies -- but since all $Q_a$ are finite 
the assertion that
$G_a$ be of type $F_n$ is equivalent to $N$ being of type $F_n$. Hence
Theorem A becomes
~\\

{\bf Theorem A$'$.} {\em Let Q act cocompactly and properly discontinuous on
  M. Then} $\rho$ {\em is} $CC^{n-1}$ {\em if and only if N has type}
$F_n$.

\hfill $\Box$
~\\

2.4 {\bf The Openness Theorem.} The main result in this first
paper is concerned with
general -- not necessarily discrete -- actions.
~\\

{\bf Theorem B.} {\em Let (M,d) be a proper} CAT(0) {\em space, let n be a
  non-negative integer and let G be a group of type} $F_n$. {\em
  The set of all isometric actions of G on M which are cocompact and}
$CC^{n-1}$ {\em is an open subset of} $\Hom(G,\Isom(M))$.
~\\

Here, $\Isom(M)$ is the topological group of isometries of $M$ and
$\Hom(G,\Isom(M))$ is the space of homomorphisms of 
the discrete group $G$ into $\Isom(M)$;
both function spaces have the compact-open topology. Theorem B is Theorem
7.7 below. The case $n = 0$ says that cocompactness is an open condition
since every $\rho$ is $CC^{-1}$ (though not uniformly:  see Proposition
3.2.).
~\\

{\bf Corollary C.} {\em Let} ${\cal R}(G,M)$ {\em denote the space of
  all isometric actions of G on M which have discrete orbits. 
Then the set of all
  isometric actions} $\rho \in {\cal R}(G,M)$ {\em which are cocompact
  and have point stabilizers of type} $F_n$ {\em is open in} ${\cal
  R}(G,M)$.
~\\

{\bf Corollary C$'$.} {\em Let} ${\cal R}_0(G,M)$ {\em denote the subspace
  of all} $\rho \in {\cal R}(G,M)$ {\em with the property that}
$\rho(G)$ {\em acts properly discontinuously on M. Then the set of all}
  $\rho \in {\cal R}_0(G,M)$ {\em which are cocompact and have} $\ker
    \rho$ {\em of type} $F_n$ {\em is open}\footnote{Corollary C$'$ has
predecessors in the literature for the case of homomorphisms $\rho : G \to
\ZZ$.  Openness of the condition ``ker $\rho$ is finitely generated'' was
proved in [Ne 79], and of the condition ``ker $\rho$ is finitely presented''
in [FrLe 85]. See also [BRe 88] and [Re 88]} {\em in} ${\cal R}_0(G,M)$. 
~\\

There is no hope of a general openness result in $\Hom(G,\Isom(M))$ for
the finiteness properties ``$\ker \rho$ is of type $F_n$'' or ``the point
stabilizers of $\rho$ are of type $F_n$''. This indicates the
advantage of the property $CC^{n-1}$  over these
traditional finiteness properties. To get a counterexample, consider a
finitely generated group $G$ whose Abelianization $G/G'$ is free of rank
2, and take $M$ to be the Euclidean line. Then every non-discrete
translation action of $G$ on $\EE^1$ has kernel the commutator subgroup
$G'$. But the non-discrete translation actions are dense in the space of
all translation actions. So if we had an openness result for the
property ``$\ker \rho$ is finitely generated'', it would imply ``$G'$ is
finitely generated if (and only if) some homomorphism $\chi: G \dpo \ZZ$
has finitely generated kernel''. This is absurd as is shown by the direct
product $G = \langle a,x|xax^{-1} = a^2\rangle \times \ZZ$ which has
commutator subgroup isomorphic to the dyadic rationals, i.e. $G' \cong \ZZ
[\frac{1}{2}]$. 
~\\

2.5 {\bf Connections with Lie groups and local rigidity.} The following
examples explain how our openness results Theorem B and Corollaries C
and C$'$ are related to locally rigid isometric actions of discrete
groups on classical symmetric spaces.

{\bf Example}: Let $M$ be a locally symmetric space of non-compact
type (e.g. the quotient of a virtually connected non-compact linear
semisimple Lie group by a maximal compact subgroup). The natural
Riemannian metric makes $M$ a proper CAT(0) space. The
group $\Isom (M)$ is a Lie group. Call its Lie algebra $\mathfrak{g}$. Each
representation $\rho \in \Hom(G,\Isom(M))$ makes $\mathfrak{g}$ into a $\ZZ
G$-module which we denote by $\mathfrak{g}(\rho)$. A theorem of Weil 
[We 64] says
that if $G$ is finitely generated and if $H^1(G;{\mathfrak g}(\rho)) = 0$ then
all nearby representations are conjugate to $\rho$ in $\Isom(M)$,
i.e. $\rho$ has a neighbourhood $N$ in $\Hom(G,\Isom(M))$ such that
every $\rho' \in N$ is of the form $\rho'(g) = \gamma\rho(g)\gamma^{-1}$
where $\gamma$ (dependent on $\rho'$) is an isometry of $M$; then $\rho$
is said to be {\em locally rigid} (see [Ra p.90]). In that case
$\ker(\rho') = \ker(\rho)$ for all $\rho' \in U$ -- a much stronger
statement than the conclusion of Corollary C$'$. But Corollary C$'$
holds in situations where $H^1(G;{\mathfrak g}(\rho)) \not= 0$, so one may
wish to think of it as a weak form of local rigidity: the kernels may
not be locally constant, but their finiteness properties are locally
constant. The next example illustrates this:

{\bf Example}: Let $G$ be the group presented by $\langle x,y | xy^2
=y^2x\rangle$. For $n \geq 0$ define $\rho_n: G \to \ZZ$ by $\rho_0(x) =
0, \rho_0(y) = 1$, and when $n \geq 1\; \rho_n(x) = n, \; \rho_n(y) =
n^2$. It is shown in [BS] that $\ker(\rho_0)$ is a free
group of rank 2 and when $n \geq 1, \ker(\rho_n)$ is a free group of
rank $n^2 + 1$. For $n \geq 1$ define $\tilde{\rho}_n: G \to \RR$ by
$\tilde{\rho}_n(g) = \frac{1}{n^2} \rho_n(g)$. Identifying $\RR$ with the
translation subgroup of $\Isom(\RR)$, we see that 
$\{\tilde{\rho}_n\}$ converges
to $\rho_0$ in $\Hom(G,\Isom(\RR))$, and 
$\ker(\tilde{\rho}_n) = \ker(\rho_n)$ for all $n \geq 1$. 
Indeed, each $\tilde{\rho}_n$ is a
cocompact action and $\tilde{\rho}_n(G)$ acts properly discontinuously on
$\RR$. This is a case where Corollary C$'$ applies but local rigidity fails.
~\\

{\bf Remark}:  The paper [Fa 99] contains results in a Lie group context
which can be seen as analogous to Theorem B and Corollary C$'$.
~\\

2.6 {\bf The new tool.} If $X$ and $Y$ are two CW complexes, we
 write $\hat{{\cal F}}(X,Y)$ for
the set of all cellular maps $f: D(f) \to Y$, where $D(f)$ is a finite
subcomplex of $X$. By a {\em sheaf of maps on X with values in Y} we 
mean\footnote{A more leisurely exposition of all this is given in \S4.
Here, we say just enough about sheaves and finitary maps to state Theorem
D.}
any subset ${\cal F}$ of $\hat{{\cal
      F}}(X,Y)$ which is closed under restrictions and finite unions. 
The sheaf ${\cal F}$ is {\em complete} (resp.
  {\em locally finite}) if each finite subcomplex of $X$ occurs as the domain
  of some member (resp. finitely many members) of ${\cal F}$. {\em
    A cross section} of the complete sheaf ${\cal F}$ is a map $X \to Y$
  whose restrictions to all finite subcomplexes lie in ${\cal F}$. Every
 cellular map $\phi: X \to Y$ is a cross section of its ``restriction'',
the sheaf $\Res(\phi)$ consisting of all restrictions of $\phi$ to finite
  subcomplexes.

These concepts become useful if $X$ and $Y$ are endowed with
cell permuting actions of a group $G$. Then $\hat{{\cal
    F}}(X,Y)$ has a natural $G$-action: If $g \in G$ and $f \in
\hat{{\cal F}}(X,Y)$ then the $g$-translate of $f$, which we write $gf
\in \hat{{\cal F}}(X,Y)$, has domain $D(gf) = g D(f)$ and maps $gx$ to $gf(x)$
for each $x \in D(f)$. A $G$-{\em sheaf} is a sheaf which is invariant under
this action. If $\phi: X \to Y$ is a $G$-equivariant cellular map then
$\Res(\phi)$ is a $G$-sheaf and is, of course, locally finite. If $\phi$
is an arbitrary cellular map then the $G$-sheaf generated by
$\Res(\phi)$ will not, in general, be locally finite. But if it is so
-- and the important fact is that this happens far beyond the
equivariant case -- we call $\phi$ a {\em finitary} (more precisely:
$G$-{\em finitary}) map. Thus a finitary map $\phi: X \to Y$ is just a
cellular map which can be exhibited as a cross section of a locally
finite $G$-sheaf.

In our situation, finitary maps will occur as cellular endomorphisms
$\phi: X^n \to X^n$ of the free, $n$-dimensional, $(n-1)$-connected,
cocompact $G$-CW complex $X^n$ of \S 2.2. Recall that $X^n$ is endowed
with a chosen $G$-equivariant control map $h: X^n \to M$ into the
CAT(0)-space $M$. The key result in the proof of Theorem B expresses the
$CC^{n-1}$ condition of $X^n$ over $a$ in terms of the following
``dynamical condition'' in $X^n$. We call a cellular map $\phi: X^n \to
X^n$ a {\em contraction (towards} $a)$ if there exists a radius $r
\geq 0$ and a number $\varepsilon > 0$ such that
\[
d(a,h \phi(x)) \leq d(a,h(x)) - \varepsilon,\; \mbox{for every}\; x \in
X^n\; \mbox{outside}\; X^n_{(a,r)}.
\]
This is independent of $a$.  We prove
~\\

{\bf Theorem D.} {\em Assume the action of G on (M,d) is
  cocompact. Then} $X^n$ {\em is} $CC^{n-1}$ {\em if and only if there
  exists a G-finitary contraction} $\phi: X^n \to X^n$.
~\\

This is contained in Theorem 6.8.
~\\

2.7 {\bf Remark on the proof of the Openness Theorem.} (This paragraph sums
up the core idea.)  The control function 
$h$ can be chosen to vary continuously with
$\rho$. Let a given cocompact action $\rho$ of $G$ on $M$ be
$CC^{n-1}$, so that we have a finitary contraction
 $\phi: X^n \to X^n$ as in Theorem
D. If we could describe $\phi$ in terms of a finite number of equations
we might expect that the very same $\phi$ would still fulfill the
properties asserted in Theorem D if the action $\rho$ were subjected to
a small perturbation. However, a description of $\phi$ requires not only
the finitary $G$-sheaf ${\cal F}(\phi)$ generated by $\Res(\phi)$ but
also an infinite number of choices of members of ${\cal F}(\phi)$.
Thus we cannot expect the same $\phi$ to work for all $\rho'$ near $\rho$.
But the sheaf ${\cal F}(\phi)$ itself can be described in terms of a
finite number of equations and we can pin down a finite
number of inequalities which are necessary and sufficient for ${\cal
  F}(\phi)$ to have a cross section which fulfills Theorem D. Thus, even
though perturbing the action $\rho$ slightly to $\rho'$ requires a new finitary
contraction $\phi'$ to establish $CC^{n-1}$ for $\rho'$, we will be able to
guarantee that $\phi'$ does exist as a cross section of the old sheaf
${\cal F}(\phi)$.
~\\

2.8 {\bf Examples.} Let $K$ be a field endowed with a discrete valuation
$v: K \to \ZZ \cup \{\infty\}$. Then we can take $M$ to be the
Bruhat-Tits-tree of $\SL_2(K)$ (see [Se], Chapitre II) acted on by
$\GL_2(K)$. Every representation $\rho: G \to \GL_2(K)$ can thus be
interpreted as an action of $G$ on $M$ with discrete orbits. 
Let $S$ be a finite set
of pairwise inequivalent discrete valuations containing $v$. Let $O_S
\subseteq K$ denote the ring of $S$-integers,\footnote{$O_S$ consists of
all $x \in K$ with $w(x) \geq 0$ for all valuations $w$ on $K$ except
possibly those in $S$.} put $G = \SL_2(O_S)$ and take $\rho_0$ to be the
natural embedding of $G$ into $\GL_2(K)$. The action $\rho_0$ is
cocompact as long as $|S|
\geq 2$.

In the case when $K$ is a finite extension of the field of rational
functions $\FF_q(t)$ over a finite field $\FF_q$ we know by a result
of U.~Stuhler [St 80] that $G$ is of type $F_{|S|-1}$ but not of type
$F_{|S|}$. This applies also to the point stabilizers of $\rho_0$ which are 
$\SL_2(O_S) \cap \GL_2(O_v) = \SL_2(O_{S-\{v\}})$, and hence the point
stabilizers of $\rho_0$ are of type
$F_{|S|-2}$ but not $F_{|S|-1}$. By Theorem A we conclude that $\rho_0$ is
$CC^{|S|-3}$ but not $CC^{|S|-2}$.

The interesting point here is that a similar phenomenon seems to
occur in the parallel case when $S$ is a finite set of rational primes
and $G = \SL_2(\ZZ_S)$ acts by Moebius transformations on the
hyperbolic plane $M = \HH^2$. This situation is of course more subtle
since this action does not have discrete orbits
when $|S| \geq 1$ so that Theorem A is
not applicable. Nevertheless we conjecture\footnote{One expects the
finiteness properties of $S$-arithmetic groups to be quite
different over function fields than over number fields. For instance,
``$F_{k-1}$ but not $F_k$'' does not occur for reductive groups in the
number field case whereas it is typical in the function field case. Our
conjecture indicates that ``$CC^{k-1}$ but not $CC^k$'' appears in the
number field case, making the two cases more similar.} that the natural
action of $G
= \SL_2(\ZZ_S)$ on $M = \HH^2$ is $CC^{|S|-2}$ but not $CC^{|S|-1}$. The
partial results we have obtained in this direction require
consideration of the $CC^{n-1}$-property over endpoints of $M$. That is
the theme of the sequel paper [BG$_{\mbox{II}}$]; this discussion of
$\SL_2(\ZZ_S)$ is continued in \S10.7(B) of that paper.

The rest of this paper consists of proofs of what has been outlined here.
An outline of the sequel paper [BG$_{\mbox{II}}$] is found in \S10.  

\newpage

\section{Technicalities Concerning the $CC^{n-1}$ Property}

3.1 {\bf The definition.} Let $(M,d)$ be a proper unique-geodesic 
space. As in \S 2.2 $X$ is a free left
$G$-complex with $G\backslash X^n$ finite and $\rho$ is a left action of $G$ on $M$ by isometries.
~\\

{\bf Proposition 3.1} {\em There exists a control function (i.e.}
$G$-{\em map} $h: X \to
M${\em ). If} $h_1$ {\em and} $h_2$ {\em are two such then}
$\sup\{d(h_1(x), h_2(x)) | x \in X^n\} < \infty$.
~\\

{\em Proof.} The map $h$ is defined inductively on skeleta. On $X^0$
define $h$ arbitrarily on one representative vertex in each orbit of
vertices, and then extend equivariantly. Assuming $h$ defined on
$X^{k-1}$, choose a representative $k$-cell $\sigma$ in each orbit of
$k$-cells. Since $M$ is contractible, for each such $\sigma$ the map
$h|\mathop \sigma\limits^\bullet$ can be extended\footnote{
$\mathop\sigma\limits^\bullet$ denotes the cell-boundary of $\sigma$,
i.e., $\mathop\sigma\limits^\bullet = \sigma \cap X^{k-1}$} to $\sigma$.
Then $h$
can be extended equivariantly to the rest of $X^k$. The second part is
clear since $G \backslash X^n$ is finite. 

\hfill $\Box$ 
~\\

Choose a control function $h: X \to M$. The property $CC^{n-1}$ was defined
in \S2.2.  
It is an immediate consequence of the triangle inequality in $M$ that if
$X$ is $CC^{n-1}$ over some $a \in M$ with lag $\lambda$ then $X$ is
$CC^{n-1}$ over any other $b \in M$  with lag $\lambda+2d(a,b)$. Hence
the property ``$CC^{n-1}$ over $a$'' is actually independent of the
choice of $a \in M$. But we have to expect that the lag $\lambda$ does
depend on $a$ (as well as on the radius $r$), and this leads us to say that
the complex $X$ is {\em uniformly} $CC^{n-1}$ {\em over} $M$ if $X$ is
$CC^{n-1}$ over every $a \in M$ with lag $\lambda$ independent of $a$.
~\\

{\bf Proposition 3.2} $X$ {\em is uniformly} $CC^{n-1}$ {\em over} $M$
{\em if and only if the given action} $\rho$ {\em of} $G$ {\em on} $M$
{\em is cocompact and} $X$ {\em is} $CC^{n-1}$ {\em over some} $a \in
M$. {\em In particular} $X$ {\em is uniformly} $CC^{-1}$ {\em if and
  only if} $\rho$ {\em is cocompact.}
~\\

{\bf Proof.} The assumption that $G\backslash X^n$ is finite applies for each $n
\geq 0$. So we can choose representatives $v_1, \ldots, v_m \in X^0$ for
the $G$-orbits, take $R$ to be the diameter of the set $h(\{v_1, \ldots,
v_m\}) \subseteq M$ and find
\[
\min d(a,G h(v_1)) \leq \min d(a,h(X^0))+R,
\]
where $d(a,S):= \{d(a,s)|s \in S\}$ for any subset $S \subseteq M$. If
$X$ is uniformly $CC^{n-1}$ over $M$ then the right hand side has a
bound independent of $a$, whence $\rho$ is cocompact. Conversely, assume
$X$ is $CC^{n-1}$ over $a$ with lag $\lambda = \lambda(a,r)$. We have 
$X_{(ga,r)} = g X_{(a,r)}$, for $g \in G$ and $r \geq 0$, so
$\lambda$ is a lag for each $a' \in
Ga$. Hence $\lambda+\min 2d(Ga,b)$ is a lag for $b \in M$, and if
$\rho$ is cocompact this has an upper bound independent of $b$. 

\hfill $\Box$
~\\

3.2 {\bf The Invariance Theorem.}
Up to now we have defined the property $CC^{n-1}$ of actions of $G$ on
$M$ by isometries using an $n$-dimensional $(n-1)$-connected cocompact
free $G$-CW-complex $X^n$ and a control function $h: X^n \to
M$. We must prove invariance: that the property is independent of the
choices of $X$ and $h$. In some
cases a natural control function presents itself on a non-free
$G$-CW-complex.\footnote{By a $G$-CW-{\it complex} we mean a CW complex 
$Y$ equipped with a 
$G$-action which permutes cells.  If this action has the additional property
that for any point $x$ in the interior of any cell $\sigma$ the stabilizers 
$G_x$ and $G_\sigma$ are equal, we call $Y$ a {\it rigid} $G$-CW-complex.}
It is useful to be able to read off the $CC^{n-1}$
property directly in such a case: 
~\\

{\bf Theorem 3.3} {\em Let G be of type} $F_n$ {\em let Y be a
  cocompact $n$-dimensional $(n-1)$-connected rigid} $G$-CW-{\em complex 
such that the
  stabilizer of each p-cell is of type} $F_{n-p},\, p \leq n-1$. {\em Let} $h:
Y \to M$ {\em be a G-map and let} $a \in M$. {\em The property that Y
  be} $CC^{n-1}$ {\em over a is independent of the choices of a, of
  Y and of the G-map h.}
~\\

In case $Y$ is a free $G$-complex, Theorem 3.3 has an elementary  
proof which we sketch below. When $Y$ is not free a more difficult proof
is required which we delay until \S9 because the methods are not related to
anything else in this paper.
~\\

{\bf Proof of Theorem 3.3 when $Y$ is free (sketch)}:
Independence of $a$ is clear. Independence of $h$ follows
from Proposition 3.1. In the free case we may always attach cells to make
the complex contractible, and extend $h$.  Let $X$ and $Y$ be two
contractible free $G$-CW complexes with cocompact $n$-skeleta.
Choose cellular maps
$\alpha: G\backslash X \to G\backslash Y$ and $\beta: G\backslash Y \to
G\backslash X$ which are mutually homotopy inverse. Their lifts
$\tilde{\alpha}: X \to Y$ and $\tilde{\beta}: Y \to X$ are bounded maps and
there is a bounded homotopy\footnote{A map [resp. homotopy] is {\em bounded}
if there is
an integer $N$ such that for all cells $\sigma$ of $\tilde X^{n-1}$, the
image of $\sigma$ [resp. of $\sigma \x I$]
lies in a subcomplex containing $\leq N$ cells.}
in $X^n$ between $\tilde{\alpha} \circ
  \tilde{\beta} | Y^{n-1}$ and the inclusion map $Y^{n-1}
  \hookrightarrow Y^n$. If $h: X \to M$ is a $G$-map then $h \circ
  \tilde{\beta}: Y \to M$ is a $G$-map. It follows directly that if $X$
  is $CC^{n-1}$ over $a$ (using $h$) then $Y$ is $CC^{n-1}$ over $a$
  (using $h \circ \tilde{\beta}$). \

\hfill $\Box$
~\\

In view of Proposition 3.2 and Theorem 3.3, the phrases ``$\rho$ is $CC^{n-1}$''
and ``$\rho$ is uniformly $CC^{n-1}$'' are unambiguous.  	 

As an application of Theorem 3.3 we find
~\\

{\bf Corollary 3.4} {\em Let} $H \leq G$ {\em be a subgroup of finite
  index in} $G$. {\em Then the G-action} $\rho$ {\em is} $CC^{n-1}$ {\em
  if and only if its restriction} $\rho|H$ {\em is} $CC^{n-1}$.
~\\

{\bf Proof.} The $G$-CW-complex $X$ and the control function 
$h: X \to M$ can also
be used to test the $CC^{n-1}$ property of the restricted action
$\rho|H$. Since the subcomplexes $X_{(a,r)}$ remain the same so do the
$CC^{n-1}$ properties. 

\hfill $\Box$
~\\

\newpage
\section{Finitary Maps and Sheaves of Maps}

This is a self-contained introduction to a new topological tool.  The
important idea---finitary maps\footnote{Seminal versions of the concepts of
finitary maps, sheaves and homotopies are contained in [BS 80], [BNS 87] and
[Re 88], specifically in the proofs of results relating the geometric
invariants to finiteness properties of normal subgroups $N\triangleleft G$
with $G/N$ Abelian.}---is introduced in \S4.7.  A finitary map is a special
kind of map between $G$-CW complexes which generalizes the notion of
equivariant map.  In our situation there are not enough equivariant maps but
there are enough finitary maps.  
~\\

4.1 {\bf Sheaves of maps.} Let $X$ and $Y$ be CW-complexes. 
By a {\em sheaf (of maps)} ${\cal F}: X \leadsto Y$ we mean a set ${\cal F}$
of cellular maps $f : D(f) \to Y$ with domain $D(f)$ a finite subcomplex of
$X$ satisfying the following axioms: 
\begin{enumerate}
\item[(i)] ${\cal F}$ contains the empty map.
\item[(ii)] If $f \in {\cal F}$ and if $K$ is a subcomplex of $D(f)$ then
$f\mid K$ is also in ${\cal F}$.  
\item[(iii)] If $f$ and $f'$ are in ${\cal F}$ and agree on the
  intersection of their domains then $f\cup f' : D(f) \cup D(f') \to Y$ is
also in ${\cal F}$.
\end{enumerate}

If $K$ is a subcomplex of $X,\; {\cal F}|K$ denotes the sheaf consisting
of all restrictions of maps in ${\cal F}$ to subcomplexes of $K$.
A {\em subsheaf} of a sheaf ${\cal F}$ is a subset which is itself a
sheaf.  Every set of maps from finite subcomplexes of $X$ to $Y$
generates a sheaf in the
obvious way. Each sheaf ${\cal F}: X \leadsto Y$ has natural
``minimal'' generators, namely the members $f \in {\cal F}$ whose
domains $D(f)$ are carriers\footnote{If $A \subset X,\; C(A)$ denotes the
  smallest subcomplex of $X$ containing $A$; it is called the {\em
    carrier} of $A$. When $A$ is compact, $C(A)$ is a finite
  subcomplex. When $\sigma$ is a cell of $X,\; C(\sigma)$ has only the
  cell $\sigma$ in the top dimension. Clearly $C(\sigma) =
  C(\mathop\sigma\limits^\bullet) \cup \sigma$ and
  $\mathop\sigma\limits^\bullet = C(\mathop\sigma\limits^\bullet) \cap
  \sigma$, so $\dim C(\mathop\sigma\limits^\bullet) < \dim C(\sigma)$.}
of cells of $X$.

An important example is the sheaf of a cellular map $\phi: X \to Y$
which we denote by $\Res(\phi)$; it consists of all restrictions of
$\phi$ to the finite subcomplexes of $X$. More generally, if $\Phi$ is
an arbitrary family of cellular maps $\phi: D(\phi) \to Y$, where each
$D(\varphi)$ is a subcomplex of $X$, we write $\Res(\Phi)$ for the sheaf
generated by the restrictions of the members of $\Phi$ to finite
subcomplexes of $X$.

By a {\em cross section of a sheaf} ${\cal F}$ we mean a cellular map
$\phi: X \to Y$ whose sheaf $\Res(\phi)$ is a subsheaf of ${\cal
  F}$. Cross sections are easy to construct when the sheaf ${\cal F}$ is
{\em homotopically closed}; by this we mean that for each $f \in {\cal F}$ and
each finite subcomplex $K \supseteq D(f)$ there is some
$\tilde{f} \in {\cal F}$ with $D(\tilde{f}) = K$ and $\tilde{f} \mid D(f) =
f$.  An easy induction on the skeleta of $K$ shows that for ${\cal 
  F}$ to be homotopically closed it suffices that for each $n$-cell
 $\sigma$ of $X$
and each $f \in {\cal F}$ with $D(f) = C(\mathop\sigma\limits^\bullet)$
there is some $\tilde{f} \in {\cal F}$ with $D(\tilde{f}) =
  C(\sigma)$ extending $f$. And with the same inductive argument
  one shows that every homotopically closed sheaf ${\cal F}:
 X \leadsto Y$ admits a cross section $\phi: X \to Y$.
~\\

4.2 {\bf $\bf G$-sheaves.} Assume that $X$ and $Y$ are
$G$-CW-complexes.  Then $G$
acts on the set of all cellular maps $\phi: K \to Y$ with $K$ a
subcomplex of $X$: if $g \in G$ we write $g\phi$ for the $g$-translate
of $\phi$; it has domain $D(g\phi) = gK$ and maps $x \in gK$ to
$g\phi(g^{-1}x)$. By a $G$-{\em sheaf} we mean a sheaf ${\cal F} : X
\leadsto Y$ which is invariant under this action. If ${\cal F}$
is a sheaf we write $G{\cal F}$ for the $G$-sheaf generated
by ${\cal F}$. Thus $G{\cal F}$ is the set of all maps 
which can be written as the union of finitely many
maps $g_if_i,\; g_i \in G,\; f_i \in {\cal F}$.

By a {\em fundamental domain} of the $G$-CW-complex $X$ we mean any
subcomplex $F \subseteq X$ with $GF = X$. A sheaf ${\cal F}_0: F \leadsto Y$
defined on a fundamental domain is $G$-{\em saturated} if whenever $f
\in {\cal F}_0$ and $g \in G$ as such that $g\, D(f) \subseteq F$
then $gf \in {\cal F}_0$. The restriction  ${\cal F}|F$ of any $G$-sheaf ${\cal
  F}: X \leadsto Y$ is certainly $G$-saturated. Conversely, $G$-saturated
sheaves lead to $G$ sheaves as follows:
~\\

{\bf Proposition 4.1.} {\em Let} $F \subseteq X$ {\em be a fundamental
  domain and} ${\cal F}_0: F \leadsto X$ {\em a sheaf. Then} ${\cal
  F}_0$ {\em is} $G$-{\em saturated if and only if} ${\cal F}_0 = G{\cal
    F}_0|F$. 
~\\

{\em Proof.} It is clear that ${\cal F}_0$ is always a subsheaf of the
$G$-saturated sheaf $G{\cal F}_0|F$. If $f \in G{\cal F}_0|F$ then $f =
g_1f_1 \cup \ldots \cup g_nf_n$ with $g_i \in G,\;
f_i \in {\cal F}_0$. Thus $g_iD(f_i) = D(g_if)
\subseteq D(f) \subseteq F$; but if ${\cal F}_0$ is $G$-saturated
this implies $g_if_i \in {\cal F}_0$ and hence $f \in {\cal
  F}_0$. 

\hfill $\Box$
~\\

Here is a useful fact about sheaves
generated by cellular maps:
~\\

{\bf Proposition 4.2.} {\em Let} $\Phi$ {\em denote a set of cellular
  maps} $\phi: D(\phi) \to Y,\; D(\phi) \subseteq X$. {\em Then we have}
$\Res(G\Phi) = G\Res(\Phi)$. {\em In particular, if} $\Phi$ {\em is}
$G$-{\em invariant, then} $\Res(\Phi)$ {\em is a} $G$-{\em
  sheaf. Moreover, in that case} $\Res(\Phi) = G\Res(\Phi|F)$ {\em for
    every fundamental domain} $F \subseteq X$.
~\\

{\em Proof.} $\Res(G\Phi)$ is generated by all restrictions $(g\phi)|K$
with $g \in G,\; \phi \in \Phi$, and $K$ a finite subcomplex of $X$. But
$(g\phi)|K = g(\phi|g^{-1}K)$, and the right hand side of this equation
exhibits generators of $G \Res(\Phi)$. This proves the first
assertion. For the second assertion note that $\Res(\Phi)$ is generated by the
restrictions of the maps $\phi \in \Phi$ to the carriers $C(\sigma)$ of
the cells $\sigma$ of $X$; and a fundamental domain $F$ will always
contain $G$-translates of these carriers. This shows that if $G\Phi =
\Phi$ then $\Res(\Phi) \subseteq G\Res(\Phi|F)$. The other inclusion is
obvious. 

\hfill $\Box$   
~\\

4.3 {\bf Locally finite sheaves.} A sheaf ${\cal F}: X \leadsto Y$
  is said to be {\em locally finite} if the restriction of ${\cal F}$ to each
  finite subcomplex $K$ of $X$ is a finite set of maps ${\cal
    F}|K: K \leadsto Y$. Note that it suffices to assume that
  the restriction of ${\cal F}$ to the carrier of each cell is finite.
~\\

{\bf Proposition 4.3.} {\em Let} ${\cal F}: X \leadsto Y$ {\em be a}
$G${\em -sheaf. If the restricted sheaf} $({\cal F}|F): F \leadsto Y$ {\em
    is locally finite for some fundamental domain} $F \subseteq X$ {\em
    then} ${\cal F}$ {\em is also locally finite.}
~\\

{\em Proof.} Let $\sigma$ be a cell of $X$, and $g \in G$ with $g\sigma
\subset F$. Since $F$ is a subcomplex it contains, in fact the whole of
the carrier $C(g\sigma) = gC(\sigma)$, and so $g({\cal F}|C(\sigma)) =
{\cal F}|C(g\sigma)$ is a subsheaf of $({\cal F}|F)$. Since ${\cal F}|F$
  is locally finite this subsheaf is finite, hence so is ${\cal
    F}|C(\sigma)$. 

\hfill $\Box$
~\\

A $G$-sheaf ${\cal F}$ is {\em finitely generated} if it is generated by
the $G$-translates of a finite subset of ${\cal F}$.

Using Proposition 4.1 one easily proves 
~\\

{\bf Proposition 4.4.} 
{\em Let $X$ and $Y$ be $G$-CW complexes and let ${\cal F} : X  
\leadsto Y$ be a $G$-sheaf.}
\begin{enumerate}
\item[a)]  {\em If ${\cal F}$ is locally finite and the $G$-action on $X$
is cocompact then ${\cal F}$ is finitely generated.}
\item[b)]  {\em If ${\cal F}$ is finitely generated and each cell of $X$ has
finite stabilizer then ${\cal F}$ is locally finite.}
\end{enumerate}

\hfill $\Box$ 
~\\

4.4 {\bf Embedding sheaves into homotopically closed sheaves.}
 In order to construct
cross sections it will often be important to embed a given sheaf in a
homotopically closed sheaf. Now, it is clear that if the complex $Y$ is contractible
then every sheaf ${\cal F}: X \leadsto Y$ can be embedded in a
homotopically closed
sheaf $\tilde{{\cal F}}: X \leadsto Y$, and if ${\cal F}$ is a $G$-sheaf
we can choose $\tilde{{\cal F}}$ to be a $G$-sheaf. Given that ${\cal
  F}$ is locally finite it requires care to make
$\tilde{{\cal F}}$ locally finite. In this paper we shall need the
embedding only in the 
situation when the $G$-$CW$-complex $X$ has finite cell stabilizers.
~\\

{\bf Proposition 4.5.} {\em Assume the} $G$-{\em complex} $X$ {\em is
  locally finite with finite cell stabilizers and} $Y$ {\em is contractible. Then every
  locally finite} $G$-{\em sheaf} ${\cal F}: X \leadsto Y$ {\em can be
  embedded in a homotopically closed locally finite} $G$-{\em sheaf} $\tilde{{\cal
      F}}: X \leadsto Y$.
~\\

{\em Proof.} By Lemma 4.6 below we can choose a fundamental domain $F
\subseteq X$ which contains only finitely
many members of each $G$-orbit of cells. Let ${\cal F}_0 = ({\cal F}|F):
F \leadsto Y$ be the restricted sheaf. ${\cal F}_0$ is locally finite
and $G$-saturated (see Proposition 4.1). If we can embed ${\cal
  F}_0$ in a locally finite $G$-saturated homotopically closed sheaf
$\tilde{{\cal F}}_0$ then $\tilde{{\cal F}} = G\tilde{{\cal F}}_0$ will
 solve our embedding problem. Indeed, by Proposition 4.1 $\tilde{{\cal F}}|F =
\tilde{{\cal F}}_0$, hence Proposition 4.3 applies to show that
$\tilde{{\cal F}}$ is locally finite; and ``homotopically closed'' is also a property
which is easily seen to be inherited from the restricted sheaf
$\tilde{{\cal F}}_0 = \tilde{{\cal F}}|F$. 

It remains to show that the locally finite and $G$-saturated sheaf
${\cal F}_0: F \leadsto Y$ can be embedded in a homotopically closed
 locally finite $G$-saturated 
sheaf $\tilde{{\cal F}}_0: F \leadsto Y$. We construct $\tilde{{\cal
        F}}_0$ inductively on the $p$-skeleton of $F$. The induction
    starts with $p = -1$ where the empty map will do. So assume $p \geq
    0$, and that a homotopically closed locally finite and $G$-saturated sheaf
    $\tilde{{\cal F}}_0^{p-1}: F^{p-1} \leadsto Y$ has been constructed
    which contains ${\cal F}_0^{p-1} = {\cal F}_0|F^{p-1}$. For each
    $p$-cell $\sigma$ of $F$ and each $f \in \tilde{{\cal
        F}}_0^{p-1}$ with domain $D(f)$ the $(p-1)$-skeleton of
    the carrier of $\sigma$, we choose a new map $\tilde{f}:
    C(\sigma) \to Y$ extending $f$ -- this is possible since $Y$
    is $p$-connected. Adjoining these maps $\tilde{f}$ to
    $\tilde{{\cal F}}_0^{p-1}$, together with all other members $f' \in
    {\cal F}^p$ with $D(f') = C(\sigma)$, yields a sheaf ${\cal G}:
    F^p \leadsto Y$ which is locally finite, homotopically closed and
 restricts to ${\cal G}|F^{p-1} = \tilde{{\cal F}}_0^{p-1}$.

${\cal G}$ is not yet $G$-saturated. In order to make it $G$-saturated
we have to adjoin, for each $p$-cell $\sigma$ of $F$, the maps $gf'$,
where $f' \in {\cal G}|C(\sigma)$ and $g$ is an element of $G$ with
$g\sigma \subseteq F$. By our careful choice of $F$ there are only
finitely many translated cells $g\sigma$ in $F$, and, since the cell
stabilizers of $X$ are finite, only finitely many elements $g \in G$
 are needed
for a given cell $\sigma$. This makes the resulting sheaf
$\tilde{{\cal F}}_0^p$ locally finite and $G$-saturated. Since
$\tilde{{\cal F}}_0^p|F^{p-1} = {\cal G}|F^{p-1}$ it is still
homotopically closed. 

\hfill $\Box$
~\\

It remains to supply the proof of
~\\

{\bf Lemma 4.6.} {\em Every locally finite} $G$-$CW${\em -complex} $X$
  {\em contains a fundamental domain} $F$ {\em with the special feature
    that} $F$ {\em contains only finitely many members of each} $G$-{\em
    orbit of cells.}
~\\

{\em Proof.} Let $T$ be a system of representatives for each $G$-orbit
of cells, and let $F \subseteq X$ be the union of all subcomplexes
$C(\sigma)$ with $\sigma$ running through $T$. Let $\tau$ be a cell in
$X$, and $g \in G$ with $g\tau \in F$. Then there is some $\sigma_g \in
T$ with $g\tau \subseteq C(\sigma_g)$, i.e., $\tau \subseteq
g^{-1}C(\sigma_g) = C(g^{-1}\sigma_g)$. But $\tau$ can only be contained
in finitely many subcomplexes of the form $C(g^{-1}\sigma_g)$. Hence
there are finitely many cells $\sigma_1, \ldots, \sigma_m \in T$ and
finitely many elements $g_1, \ldots, g_m$ with the property that for
each $g \in G$ with $g\tau \in F$, there is $j \in \{1, \ldots, m\}$
such that $g^{-1}\sigma_g = g_i^{-1}\sigma_i$ -- in other words,
$\sigma_g = \sigma_i$ and $gg_i^{-1} \in \mbox{Stab}(\sigma_i)$. It
follows that $g\tau \subseteq C(\sigma_g) = C(\sigma_i)$ for some $i$;
hence there are only finitely many possibilities for $g\tau$ to be in
$F$.

\hfill $\Box$
~\\

4.5 {\bf Composing sheaves.} In Part II we will compose sheaves
${\cal F}: X \leadsto Y, {\cal F}': Y \leadsto Z$. Now, the set ${\cal
  P}$ consisting of all compositions $f' \circ f$ with $f
  \in {\cal F},\; f' \in {\cal F}'$ and $f(D(f)) \subseteq D(f')$
  is not, in general, a sheaf -- though it does have properties (i) and
  (ii) of the definition of a sheaf in \S 4.1. We define ${\cal F}'
  \circ {\cal F}$ to be the sheaf generated by ${\cal P}$. Thus ${\cal
    F}' \circ {\cal F}$ consist of all maps $f'': D(f'') \to Z$
  which can be written as a union $f'' = f'_1f_1 \cup
  \ldots \cup f'_k f_k$ with $f_i \in {\cal F}, f'_i \in
  {\cal F}'$ and $\mbox{image} f_i \subseteq D(f'_i)$ for all $i$.
~\\

{\bf Proposition 4.7} (i) {\em Let} ${\cal F}: X \to Y$ {\em and} ${\cal G}:
Y \to Z$ {\em be sheaves. If} ${\cal F}$ {\em and} ${\cal G}$ {\em are}
$G$-{\em sheaves, so is} ${\cal G} \circ {\cal F}$. {\em If} ${\cal F}$
 {\em and} ${\cal G}$ {\em are locally finite, so is} ${\cal G} \circ {\cal
    F}$.\\
(ii) {\em If} $\phi: D(\phi) \to Y$ {\em and} $\psi: D(\psi) \to Z$ {\em
  are cellular maps defined on subcomplexes} $D(\phi) \subseteq X,
D(\psi) \subseteq Y$ {\em with} $\phi(D(\phi)) \subseteq D(\psi)$. {\em
  Then} $\Res(\psi \circ \phi) = \Res(\psi) \circ \Res(\phi)$.
~\\

Both assertions are straightforward. 

\hfill $\Box$
~\\

4.6 {\bf Homotopy of sheaves.} Let $X$ and $Y$ be CW-complexes. We consider the
unit interval $I$ and the product $X \times I$ with their canonical
CW-structures. A {\em homotopy of sheaves}, ${\cal H}: X \times I
\leadsto Y$ is a set of cellular maps $H: D(H) \to Y$, where $D(H) = D_H
\times I$ with $D_H$ a finite subcomplex of $X$, such that ${\cal H}$
contains the empty maps and is closed under union and
restriction (to subcomplexes of the form $D \times I$). More precisely,
we will say that such a set of maps ${\cal H}$ is a {\em homotopy between the
sheaves} ${\cal F}_0 = {\cal H}|X \times \{0\}$ and ${\cal F}_1 = {\cal
  H}|X \times \{1\}$, where ${\cal F}_0$ and ${\cal F}_1$ are regarded as
sheaves $X \leadsto Y$.

Although a homotopy of sheaves ${\cal H}$ is not a sheaf in
the technical sense (restrictions to subcomplexes other than those of
the form $D \times I$ are not considered) we can still use constructions
like unions, restrictions and compositions, and we can talk about
$G$-homotopies in the obvious sense.

An obvious necessary condition for the existence of a homotopy between
two sheaves ${\cal F}_0, {\cal F}_1: X \leadsto Y$ is that the collections
of domains of ${\cal F}_0$ and of ${\cal F}_1$ coincide. If they do,
then contractibility of $Y$ is sufficient for the existence of a
homotopy ${\cal H}: {\cal F}_0 \simeq {\cal F}_1$. 
Parallel to Proposition 4.5 we find
~\\

{\bf Proposition 4.8.} {\em Assume the} $G$-{\em complex} $X$ {\em is
  locally finite with finite cell stabilizers and} $Y$ {\em is
 contractible. Then any two
  locally finite} $G$-{\em sheaves} ${\cal F}_0, {\cal F}_1: X \leadsto Y$
{\em with the same collections of domains are homotopic via a locally finite}
$G$-{\em homotopy} ${\cal H}: X \times I \leadsto Y$.

\hfill $\Box$
~\\

4.7 {\bf Finitary maps.} A cellular map $\phi: X \to Y$ between
two $G$-CW-complexes is said to be {\em finitary} (or, more precisely,
$G$-{\em finitary}) if $\phi$ is a cross section of a locally
finite $G$-sheaf ${\cal F}: X \leadsto Y$. Such a $G$-sheaf ${\cal F}$
will have to contain all of $G\Res(\phi)$; hence $G\Res(\phi) =
\Res(G\phi)$ is locally finite. This shows that $\phi: X \to Y$ is
finitary if and only if the $G$-sheaf $G \Res(\phi) = \Res(G\phi)$ is
locally finite.
Hence, in view of Proposition 4.7 we obtain
~\\

{\bf Proposition 4.9.} {\em If} $\phi: X \to Y$ {\em and} $\psi: Y \to Z$
{\em are}  $G$-{\em finitary maps, so is their composition} $\psi \circ
\phi: X \to Z$. 

\hfill $\Box$
~\\

Examples of finitary maps are, of course, the $G$-equivariant maps
$\phi: X \to Y$; for if $\phi$ is equivariant then $\Res(\phi)$ is a
$G$-sheaf with cross section $\phi$. But there are many finitary
maps beyond the equivariant ones. Take, for instance, $X$ and $Y$ to be
the 0-dimensional free $G$-space $G$. A finitary map $\phi: G \to G$ is
then given by a finite set of elements $T \subseteq G$ and a function
$\kappa: G \to T;\; \phi(g) = g\kappa(g),\; g \in G$. If one wishes to
extend this map $\phi$ to a finitary map $\phi: \Gamma \to \Gamma$ on 
the Cayley graph of $G$
with respect to a set of generators $S \to G$, one has to choose a
finite set of edge paths $P$ with origin and terminus in $T$ and a map
$\mu: G \times S \to P$ such that $\mu(g,s)$ has origin $\kappa(g)$
and terminus $\kappa(gs)$, and put $\phi((g,s)): = \mu(g,s)$.
~\\

{\bf Proposition 4.10.} {\em If the} $G$-$CW$-{\em complex} $X$ {\em is
  locally finite with finite cell stabilizers and} $Y$ {\em is
 contractible then any two}
$G${\em -finitary maps} $\phi_0,\; \phi_1: X \to Y$ {\em are homotopic via
    a} $G$-{\em finitary homotopy.}
~\\

{\em Proof.} Take locally finite sheaves ${\cal F}_0, {\cal F}_1: X
\leadsto Y$ with cross sections $\phi_0, \phi_1$, respectively. By
Proposition 4.8 there is a locally finite homotopy ${\cal H}: X \times I
\leadsto Y$ between ${\cal F}_0$ and ${\cal F}_1$; ${\cal H}$ is a
$G$-homotopy. By an inductive construction one finds a homotopy $\psi: X
\times I \to Y$ from $\phi_0$ to $\phi_1$ which is a cross section of
${\cal H}$. The details are left to the reader. 

\hfill $\Box$
~\\

\newpage
\section{Sheaves and Finitary Maps Over a Control Space}

In this section we consider sheaves ${\cal F}: X \leadsto X$ of maps on
a CW-complex into itself over a ``control space'' $M$. Throughout the
section $M$ is a proper metric space -- the CAT(0) property will not be
used. The group $G$ is assumed to act by cell permuting automorphisms on
$X$ and by isometries on $M$, and we are given a control
function $h: X \to M$.
~\\

5.1 {\bf Displacement function and norm.} Given a map $f: D(f) \to
X$ with $D(f)$ a (not necessarily finite) subcomplex of $X$ we consider
\begin{itemize}
\item the {\em displacement function of} $f$ {\em over} $M,\; \alpha_f:
  D(f) \to [0,\infty)$
defined by $\alpha_f(x): = d(h(x), hf(x)),\; x \in D(f)$, and
\item the {\em norm of} $f$ {\em over} $M,\; \|f\| \in [0,\infty]$, defined by $\|f\|: = \sup \alpha_f(D(f))$.
\end{itemize}

The displacement function of $f$ is a continuous non-negative
function. As $G$ acts by isometries and $h$ is a $G$-map it satisfies
\[
(5.1) \qquad \alpha_{gf}(gx) = \alpha_f(x), \qquad \mbox{all}\quad g \in G,\; x \in
X.
\]
The norm $\|f\|$ may be infinite if $D(f)$ is not finite. If $\|f\| <
\infty$ we call $f$ {\em bounded} (over $M$). We put
$\|\emptyset\| = 0$ if $\emptyset$ is the empty map.

It is an easy matter to verify that the following formulae hold for all
maps $f_i: D(f_i) \to X$ as above and all $g \in G$ 
\begin{enumerate}
\item[(5.2)] $\|g\; f_i\| = \|f_i\|$,
\item[(5.3)] $\|f_1 \cup f_2\| = \max(\|f_1\|, \|f_2\|)$, provided $f_1
  \cup f_2$ exists,
\item[(5.4)] $\|f_1 \circ f_2\| \leq \|f_1\| + \|f_2\|$.
\end{enumerate}

The notion of displacement and norm extend to locally finite sheaves. If
${\cal F}: X \leadsto X$ is a locally finite sheaf and $D({\cal F})$ 
stands for the union of all domains $D(f),\; f \in {\cal F}$, then
$\alpha_{{\cal F}}: D({\cal F}) \to [0,\infty)$ is the continuous function
 $\alpha_{\cal F}(x) = \sup\{\alpha_f(x)|f \in {\cal F}\}$ and
$\|{\cal F}\| :=
 \sup \alpha_{\cal  F}(D({\cal F}))$. Again we call ${\cal F}$ {\em bounded}
 if $\|{\cal F}\| < \infty$.

(5.3) and (5.4) show that if ${\cal F}$ and ${\cal G}$ are
locally finite sheaves $X \leadsto X$ then the norm of their composition
satisfies
\[
(5.5)\qquad \|{\cal F} \circ {\cal G}\| \leq \|{\cal F}\| + \|{\cal
  G}\|.
\]
Useful elementary facts concerning $G$-sheaves $X \leadsto X$ and finitary
maps $\phi : X \to X$ are collected in
~\\

{\bf Proposition 5.1} {\em (a) If} ${\cal F}: X \leadsto X$ {\em is a
  G-sheaf then the displacement function} $\alpha_{\cal F}$ {\em satisfies} $\alpha_{\cal
    F}(gx) = \alpha_{\cal F}(x)$ {\em for every} $x \in D({\cal F})$.\\
\qquad {\em (b) If} $G\backslash X$ {\em is compact then every locally
  finite G-sheaf ${\cal F}$ is bounded, and if $\phi : X \to X$ is a cross
section of ${\cal F}$ then $||\phi|| \leq ||{\cal F}||$.}
~\\

{\em Proof.} (a) is immediate from (5.1); (5.2) and
(5.3) show that the norm of a $G$-sheaf is given by the norm of a
$G$-generating subset, hence (b). 

\hfill $\Box$
~\\

5.2 {\bf Shift towards a point of $\bf M$.} Let $a \in M$.  The 
  {\em shift} of a map $f: D(f) \to X,\; D(f)$ a subcomplex of $X$,
is defined to be the continuous function $\sh_{f,a}: D(f) \to \RR$ given
by
\[
\sh_{f,a}(x): = d(h(x),a) - d(hf(x),a), \quad x \in D(f).
\] By the triangle inequality we have
\[
(5.6) \qquad |\sh_{f,a}(x)\; |\; \leq \alpha_f(x), \quad \mbox{all} \quad x
\in D(f).
\]
Hence the shift function has the global bound $\|f\|$ if $f$ is a bounded
map.

If $a,b$ are two points of $M$ then their shift functions are
related by the inequality
\[
(5.7) \qquad |\sh_{f,a}(x) - \sh_{f,b}(x)\; |\; \leq 2d(a,b).
\]
~\\

5.3 {\bf Contractions.} Let $a \in M$. We call a cellular map
 $\phi: X \to X$ a
{\em contraction of} $X$ {\em towards} $a$ if there are numbers $R \geq 0,\;
\varepsilon > 
0$ with the property that $\sh_{\phi,a}(x) \geq \varepsilon$ for every
$x \in X$ with $d(h(x),a) \geq R$. Any such numbers $R,\; \varepsilon$
are called an {\em event radius} and an {\em almost guaranteed shift} for
$\phi$. Using (5.7) one observes that, for each $b \in M$ with $d(a,b) <
\frac{\varepsilon}{2}, \quad \phi$ will then also be a contraction
towards $b$ (with event radius $R + d(a,b)$ and almost guaranteed shift
$\varepsilon - 2d(a,b)$).
~\\

{\bf Proposition 5.2} {\em If} $X$ {\em admits a contraction towards} $a
\in M$ {\em with event radius} $R$ {\em then} $X$ {\em is} $CC^{-1}$ {\em
  over} $a$ {\em with lag} $R$.
~\\

{\em Proof.} Let $\varepsilon$ be an almost guaranteed shift of this
 contraction
$\phi: X \to X$. Pick any $x \in X^0$ and $m \in \NN$ with $m \cdot
\varepsilon \geq d(a,h(x)) - R$. Then one of $h(x), h \phi(x), \ldots, h
\phi^m(x)$ is in $B_R(a)$, hence the assertion.

\hfill $\Box$
~\\

{\bf Proposition 5.3} {\em Let} $\phi: X \to X$ {\em be a bounded
  contraction towards} $a \in M$ {\em with event radius} $R$ {\em and
  almost guaranteed shift} 
$\varepsilon > 0$. {\em Then}
\[
d(h \phi^m(x),a) \leq \max(d(h(x),a) - m \varepsilon,\; R + \|\phi\|)
\]
{\em for each} $m \in \NN$ {\em and} $x \in X$. {\em In particular}
$\phi^m: X \to X$ {\em is a contraction with event radius} $R + \| \phi
\| + m \varepsilon$ {\em and almost guaranteed shift} $m\varepsilon$.
~\\

{\em Proof.} Induction on $m$ using the fact that $d(h \phi(y),a) =
d(h(y),a)-\sh_{\phi,a}(y)$ is at most $d(h(y),a)-\varepsilon$, if $d(h(y),a)
\geq R$, and at most $R + \|\phi\|$ if $d(h(y),a) \leq R$. 
 
\hfill $\Box$
~\\

{\bf Corollary 5.4} {\em If} $a$ {\em contraction} $\phi: X
\to X$ {\em towards} $a \in M$ {\em exists, then there are
  contractions} $\phi_b: X \to X$ {\em towards each} $b \in M$. {\em If} $G$
 {\em acts cocompactly on} $M$ {\em then these}
$\phi_b$ {\em can be chosen with uniform event radius} $R$. {\em If}
$\phi$ {\em is finitary these} $\phi_b$ {\em can be chosen to be finitary.}
~\\

{\bf Proposition 5.5} {\em Assume that the G}-CW-{\em complex} $X$ {\em is
  contractible, with finite cell stabilizers and cocompact} $n$-{\em
  skeleton. If} $X^n$ {\em admits a finitary
  contraction} $\phi: X^n \to X^n$ {\em towards} $a$
{\em with event radius} $R$ {\em and
  almost guaranteed shift} $\varepsilon > 0$ {\em then there is a
  cellular deformation} $\psi: X^n \times [0,\infty) \to X^{n+1}$ {\em and a
  number} $\lambda \geq 0$ {\em such that}
\[
d(h\Psi(x,s),a) \leq \max\Big(d(h(x),a) - s\varepsilon,\; R + \|\phi\|\Big)
+ \lambda
\]
{\em for every} $x \in X^n$ {\em and} $s \in [0,\infty)$.
~\\

{\em Proof.} By Proposition 4.10 there is a $G$-finitary homotopy
$\psi_0: X^n \times I \to X^{n+1}$ between  $\Id_{X^n}$ and $\phi$.
By analogy with the case of maps, let $||\psi_0|| := \ds{\sup_{t\in
I}}||\psi_0(\cdot,t)||$.  Since $\psi_0$ is finitary, $||\psi_0|| < \infty$.

Write $s \in [0,\infty)$ as $s = m+t$ where $m$ is a non-negative
integer and $t \in [0,1)$. Define $\psi: X \times [0,\infty) \to X$
to be the (continuous) cellular map $\psi(x,s) =
\psi_0(\phi^m(x),t);\quad y:= (\phi^m(x),t)$ and $z:= (\phi^m(x),0)$ are
two points of $X \times I$.  Hence
\[
\begin{array}{lcl}
d(h\psi_0(y),\; h\psi_0(z)) & \leq & d(h\psi_0(y),\; h(\phi^m(x))) + 
d(h(\phi^m(x)),\;
\psi_0(z))\\
& \leq & 2\|\psi_0\|.
\end{array}
\]
It follows by using Proposition 5.3 that
\[
\begin{array}{lcl}
d(h\psi(x,s),a) & \leq & d(h\phi^m(x),a) + 2\|\psi_0\|\\
& \leq & \max(d(h(x),a) - m\varepsilon,\; R + \|\phi\|) + 2\|\psi_0\|\\
& \leq & \max(d(h(x),a) - s\varepsilon,\; R + \|\phi\|) +
2\|\psi_0\|+\varepsilon,
\end{array}
\]
so we can choose $\lambda$ to be $2\|\psi_0\|+\varepsilon$.

\hfill$\Box$
~\\

{\bf Theorem 5.6} {\em If} $X$ {\em is contractible with finite cell
  stabilizers  and cocompact} $n$-{\em skeleton then the existence of a
  finitary contraction} $\phi: X^n \to X^n$ {\em towards} $a \in M$ 
{\em implies that} $X$ {\em is} $CC^{n-1}$ {\em over} $a$ {\em with 
constant lag}.
~\\

{\em Proof.} Let $f: S^p \to X^p_{(a,t)}$ be a map, $0 \leq p \leq
n-1$. Since $X^{p+1}$ is $p$-connected $f$ extends to a map $f_1:
B^{p+1} \to X^{p+1}$. Applying Proposition 5.3, the map $\psi \circ (f
\times \Id): S^p \times [0,\infty) \to X^{p+1}$ can be ``glued'' to
$f_1$ across $S^p \times \{0\}$ to give a map $f_2: \RR^{p+1} \to
X^{p+1}$. Write $\RR^{p+1} \times [0,\infty) = \RR_+^{p+2}$
and consider $\psi \circ (f_2 \times \id): \RR_+^{p+2} \to X^{p+2}$. The
annulus in $\RR^{p+1}$ with inner radius 1 and outer radius $r$ (sufficiently
large) adjoined along the sphere of radius $r$ to the hemisphere of
radius $r$ in $\RR^{p+1}$ gives a $(p+1)$-ball in $\RR_+^{p+2}$ whose
boundary is $S^p \subseteq \RR^{p+1} \subseteq \RR_+^{p+2}$. The
restriction of $\psi \circ(f_2 \times \id)$ to this $(p+1)$-ball
$\tilde{B}^{p+1}$ is the required $\tilde{f}: \tilde{B}^{p+1} \to
X^{p+1}$ of $f$. See Figure I. The lag is as $\lambda$ in Proposition
5.5.

\hfill $\Box$
~\\
\vspace{10cm}

\centerline{Fig. I}
~\\

5.4 {\bf Guaranteed shift.} By the {\em guaranteed shift towards} $a \in
M$ of a bounded cellular map $f: D(f) \to X$ on a subcomplex $D(f)$ of
$X$ we mean the real number
\[
\gsh_a(f): = \inf \sh_{f,a} (D(f)).
\]
Observe that $\gsh_a(f)$ is compatible with the $G$-action in the sense
that 
\[
\gsh_{ga}(gf) = \gsh_a(f)
\]
for each $g \in G$.  By (5.7) $\gsh_a(f)$ is continuous in the variable $a$
when $D(f)$ is compact.

Let ${\cal F}: X \leadsto X$ be a locally finite homotopically 
closed $G$-sheaf. The
assumption that ${\cal F}$ be homotopically closed guarantees that
 for every finite
subcomplex $K \subseteq X$ there are maps $f \in {\cal F}$ with $D(f)
= K$ so that we can consider the {\em maximal guaranteed shift of}
${\cal F}$ {\em on} $K$ (towards $a \in M$), defined by
\[
\mu_a ({\cal F} | K):= \max \{\gsh_a(f) | f \in {\cal F}, D(f) = K\}.
\] 
The case when $K$ is the zero skeleton of the carrier of a cell $\sigma$
of $X, K = C(\sigma)^0$, will be particularly useful to us. We call
$\mu_a({\cal F}|C(\sigma)^0)$ the {\em maximal guaranteed vertex shift
  on} $\sigma$. Abusing notation we shorten this to $\mu_a({\cal
  F}|\sigma)$. As to the $G$-action we have compatibility in the sense that
\[
\mu_{ga}({\cal F}|g\sigma) = \mu_a({\cal F}|\sigma)
\]
for all $g \in G, a \in M$.  As above, $\mu_a({\cal F}\mid \sigma)$ is
continuous in $a$.  
~\\

5.5 {\bf Defect of a sheaf.} Again let ${\cal F}: X \leadsto X$ be a
locally finite homotopically closed $G$-sheaf. Let $\sigma$ be
a cell of $X$, with $\dim\sigma \geq 1$. We wish to compare the
guaranteed shift $\gsh_a(f)$ of a map $f \in {\cal F}$ whose domain is
the carrier $C(\mathop\sigma\limits^\bullet)$ with the guaranteed
 shift of the extensions
$\tilde{f}: C(\sigma) \to X$ of $f$ which lie in ${\cal F}$: We define the
{\em defect} of ${\cal F}$ (towards $a$) on $\sigma$ to be the number
~\\

\noindent
$d_a({\cal F}|\sigma) = \max\limits_{f\in{\cal F}}\;
\min\limits_{\tilde{f}\in{\cal F}} \{ \gsh_a(f) - \gsh_a(\tilde{f})|D(f)
= C(\mathop\sigma\limits^\bullet), D(\tilde{f}) = C(\sigma), \tilde{f}$
extends $f\}$.
~\\

Here again $d_a({\cal F}|\sigma)$ is a shortening of what should be
$d_a({\cal F}|C(\sigma))$.

Since $\tilde f$ extends $f$ this is a non-negative number. To say that
$d_a({\cal F}|\sigma) \leq k$ is to say that for every $f \in {\cal F}$
with $D(f) = C(\mathop\sigma\limits^\bullet)$ there is an extension
$\tilde{f}$ 
of $f$ in ${\cal F}$ with $D(\tilde{f}) = C(\sigma)$ such that
$\gsh_a(\tilde{f}) \geq \gsh_a(f)-k$. Thus $k$ is an upper bound for the
``loss of guaranteed shift'' towards $a$ in extending $f$ from
$C(\mathop\sigma\limits^\bullet)$ to $C(\sigma)$ using maps in ${\cal F}$;
$d_a({\cal F}|\sigma)$ is the best upper bound.

Again, we observe that $d_a({\cal F}|\sigma)$ is compatible with the
$G$-action in the sense that 
\[
d_{ga}({\cal F}|g\sigma) = d_a({\cal F}|\sigma)
\]
for every $g \in G, a \in M$ and all cells $\sigma$ of $X$.  As above,
$d_a({\cal F}\mid \sigma)$ is continuous in $a$. 

Guaranteed vertex shift and defect are used to control the shift towards
$a \in M$ of a cross section $\phi: X \to X$ of ${\cal F}$, which we 
construct skeleton by
skeleton as follows. For each vertex $v \in X^0$ we choose some $f \in
{\cal F}$ with $v \in D(f)$ and $\sh_{f,a}(v) = {\rm gsh}_a(f|\{v\}) =
\mu_a({\cal F}|\{v\})$, and put $\phi_0(v):= f(v)$. This constructs
$\phi_0: X^0 \to X^0$ with the property that $\gsh_a(\phi_0|K^0) =
\mu_a({\cal F}|K^0)$ for every finite subset $K^0 \subseteq X^0$. Assume
now a cross section $\phi_{k-1}$ of ${\cal F}|X^{k-1}$ has been 
constructed for some $k \geq 1$. If $\sigma$ is a $k$-cell of $X$ we can
extend $\phi_{k-1}$ to the carrier $C (\sigma)$ by a map $\tilde{f}:
C(\sigma) \to X$ in ${\cal F}$ with $\gsh_a(\tilde{f}) \geq
\gsh_a(\phi_{k-1}|C(\mathop \sigma\limits^\bullet)) - d_a({\cal
  F}|\sigma)$. This constructs a cross section $\phi_k$ of ${\cal
  F}|X^k$ with the property that for each finite subcomplex $K \subseteq
X^k$
\[
\gsh_a(\phi_k|K) \geq \gsh(\phi_{k-1}|K \cap X^{k-1}) - \max\{d_a({\cal
    F}|\sigma)|\sigma\; \mbox{a}\; k\mbox{-cell of}\; K\},
\]
hence,
\[
\gsh_a(\phi_k|K) \geq \mu_a({\cal F}|K^0) - \sum\limits^k_{j=1} \max
\{d_a({\cal F}|\sigma)|\sigma\;\mbox{a}\; j\mbox{-cell of}\; K\}.
\]
For every cell $\sigma$ of
$X$ we define the {\em total defect of} ${\cal F}$ {\em on} $\sigma$
(towards $a \in M$) to be
\[
\delta_a({\cal F}|\sigma):= \sum\limits ^{\dim\sigma}_{j=1}
\max\{d_a({\cal F}|\tau)|\tau\; \mbox{a}\; j\mbox{-cell of}\;
 C(\sigma)\}.
\]
If $\dim \sigma = 0$ this means that $\delta_a({\cal F}|\sigma) = 0$. We
have proved
~\\

{\bf Proposition 5.7.} {\em Let} $a \in M$. {\em Then every locally
  finite homotopically closed sheaf} ${\cal F}: X \leadsto X$ {\em has a cross
  section} $\phi: X \to X$ {\em with}
\[
\sh_{\phi,a}(x) \geq \mu_a({\cal F}|\sigma) - \delta_a({\cal
  F}|\sigma),
\]
{\em for each cell} $\sigma$ {\em of} $X$ {\em and each} $x \in \sigma$.

\hfill $\Box$ 
~\\

{\bf Remark 5.8} In the construction proving Proposition 5.7 we have not
used the $G$-action. But if ${\cal F}$ is, in fact, a $G$-sheaf (so that
$\phi$ is finitary) and the
stabilizer of $a \in M$, denoted by $G_a$, acts freely on $X$ then in each
step of the construction the maps $\phi_k$ and hence also the cross
section $\phi$ in Proposition 5.7 can be chosen to be $G_a$-maps. 

To make Proposition 5.7 useful we must arrange for $\mu_a({\cal
  F}|\sigma)-\delta_a({\cal F}|\sigma)$ to be positive.

\newpage
\section{Construction of Sheaves with Positive Shift}

In this section $(M,d)$ is a proper unique-geodesic metric $G$-space (the
need for CAT(0) first appears in \S6.3) and $X$ is a
$G$-CW-complex with a control function $h: X \to M$.
 In view of Proposition 5.7. we wish to
construct locally finite closed $G$-sheaves ${\cal F}: X \leadsto X$
with the property that there are numbers $R \geq 0, \varepsilon > 0$ with
\[
\mu_a({\cal F} | \sigma) - \delta_a({\cal F} | \sigma) \geq \varepsilon
\]
for every cell $\sigma$ of $X$ with $h(\sigma) \cap B_R(a) =
\emptyset$. Proposition 5.7. will then show that the sheaf ${\cal F}$
admits a cross section $\phi: X \to X$ which is a contraction in the sense
of \S5.3.
~\\

6.1. {\bf The case when dim X = 0}.
~\\

{\bf Theorem 6.1} {\em Assume} $X = X^0$ {\em is a discrete} $G$-{\em
  set with finite stabilizers and finitely many orbits. Then the
  following conditions are equivalent}.
\begin{enumerate}
\item[(i)] {\em The} $G$-{\em action on M is cocompact}\footnote{Recall
    that by Proposition 3.2 this is equivalent to uniformly $CC^{-1}$.}
\item[(ii)] {\em For every} $a \in M$ {\em there is a G-finitary
    contraction} $\phi_a: X \to X$ {\em with event radius} $R$ {\em
    independent of a.}
\item[(iii)] {\em There are numbers} $R \geq 0$ {\em and} $\alpha > 0$
  {\em with the property that for each} $a \in M$ {\em one can construct
    a locally finite G-sheaf} ${\cal F}: X \leadsto X$ {\em satisfying}
    $\mu_a({\cal F}|\{x\}) \geq \alpha$ {\em for all} $x \in X$ {\em
      with} $d(h(x),a) \geq R$.
\item[(iv)] {\em For any given} $\alpha > 0$ {\em there is some} $R =
    R(\alpha)$ {\em such that (iii) holds.}
\end{enumerate}
~\\

{\em Proof.} The implication (ii) $\Rightarrow$ (i) is covered by
Proposition 5.2. The implication (iii) $\Rightarrow$ (ii) is
covered by Proposition 5.7. As (iv) $\Rightarrow$ (iii) is trivial it
remains to prove (i) $\Rightarrow$ (iv).

Let $F \subseteq X$ be a fundamental domain and $s \geq 0$
a radius with the property that $GB_s(b) = M$ for every choice of $b \in
M$. By replacing $F$ by a $G$-translate and increasing $s$ if necessary
we may assume $a \in h(F) \subseteq B_s(a)$.

Let $R = 3s+\alpha$, and let $S_{R-S}(a)$ denote the sphere of radius $R-s$
around $a \in M$. For each $u \in S_{R-S}(a)$ we choose a point $x_u \in X$
with $d(u,h(x_u)) < s$; this is possible since each point $m \in M$ has
distance at most $s$ from some point of $h(X)$. Since $d(u,h(x_u)) < s$,
$d(v,h(x_u)) < s$ for all points $v$ in a small
neighbourhood of $u$. By
\newpage 

\noindent
compactness of $S_R(a)$ we can thus arrange
that all elements $x_u$ are contained in finite subset $W \subseteq X$.

As $F$ and all stabilizers of elements of $F$ are finite, $H:= \{g \in G
| gF \cap F \not= \emptyset\}$ is a finite subset of $G$. Now, let
${\cal F}_0: F \leadsto X$ be the sheaf generated by all maps $f: \{v\}
\to X$ with $v \in F$ and $f(v) \in HW$. ${\cal F}_0$ is finite and
$G$-saturated (in the sense of \S4.2). Hence, by Propositions 4.1 and
4.3, the $G$-sheaf ${\cal F} := G{\cal F}_0$ is locally finite. If $x \in
X$ and $d(a,h(x)) \geq R$ then we pick $g \in G$ with $gx \in F$ and
consider the geodesic segment $L$ from $a$ to $ga$. One checks that
$d(ga,a) \geq R-s$; hence $L$
intersects $S_{R-S}(a)$ in a unique point $u$. Let $f(gx): = x_u$. This
defines a member $f \in {\cal F}_0$, and its translate $g^{-1}f$ takes
$x$ to $g^{-1}x_u$. We have
\[
\begin{array}{lcl}
sh_{g^{-1}f,a}(x) & = & d(h(x),a) - d(h(g^{-1}x_u),a)\\
& = & d(gh(x),ga) - d(h(x_u),ga)\\
& \geq & (d(a,ga)-s) - (d(u,ga)+s)\\
& = & d(a,ga)-s-(d(a,ga)-d(u,a))-s\\
& = & d(u,a)-2s\\
& = & R-3s = \alpha.
\end{array}
\]
This shows that $\mu_a({\cal F} | \{x\}) \geq \alpha$.

\hfill $\Box$
~\\

6.2. {\bf Measuring the loss of guaranteed shift in an extension.} 
In this subsection we
assume $n \geq 1$ and that $X$ be $CC^{n-1}$ over some point of
$M$. Then we consider an $n$-cell $\sigma$ of $X$ and a cellular map $f:
C(\mathop\sigma\limits^\bullet) \to X$ and we try to extend $f$ to a map
$\tilde{f}: C(\sigma) \to X$ with some control on the shift
$\gsh_b(\tilde{f})$ for $b \in M$.
~\\

{\bf Proposition 6.2.} {\em If} $X$ {\em is} $CC^{n-1}$ {\em over} $b \in
M$ {\em then there exists} $\lambda_0 \geq 0$ {\em and} 
$\tilde{f}: C(\sigma) \to X$
{\em such that} $\gsh_b(\tilde{f}) \geq \gsh_b(f) -
\diam hC(\sigma) - \lambda_0$.
~\\

{\bf Proof.} Let $u \in C(\mathop\sigma\limits^\bullet)$ be a point with the
property that $d(hf(u),b) =: r_0$ is the maximum distance of any point of
$hf C(\mathop\sigma\limits^\bullet)$ from $b$. The $CC^{n-1}$ hypothesis
gives a lag $\lambda(b,r)$.  Let $\lambda_0 := \lambda(b,r_0)$.   There
exists $\tilde{f}: C(\sigma) \to X$
with $h\tilde{f}(C(\sigma)) \subseteq B_{r+\lambda_0}(b)$. We then have
for each $x \in C(\sigma)$
\[
\begin{array}{lcl}
d(h(x),b) - d(h\tilde{f}(x),b) & \geq & d(h(u),b)-
\diam hC(\sigma)-d(hf(u),b)-\lambda_0\\
& \geq & \gsh_b(f) - \diam hC(\sigma)- \lambda_0,
\end{array}
\]
hence the result.

\hfill $\Box$
~\\

The same map $\tilde{f}$ works for points $b'$ near $b$ using
continuity of $\gsh_b$ in the variable $b$: given any $\varepsilon > 0$
there exist $\delta > 0$ such that
\[
\gsh_{b'}(\tilde{f}) > \gsh_{b'}(f) - \diam\, hC(\sigma) - \lambda_0 -
\varepsilon 
\]
for all $b' \in B_\delta(b)$.

We can now let $b$ range over a compact subset $B$ of $M$. Since $B$ has
finite diameter, there is a
lag $\lambda = \lambda(r)$ such that $X$ is $CC^{n-1}$ over every 
point of $B$ with
lag $\lambda$. A standard compactness argument (pass to a finite subcover
of the cover $\{B_\delta(b)\}$ of $B$) then shows
~\\

{\bf Corollary 6.3.} {\em If} $B \subseteq M$ {\em is a compact subset
  such that} $X$ {\em is} $CC^{n-1}$ {\em over
  each} $b \in B$ {\em then, given any} $\varepsilon > 0$, {\em a finite
  set} $S$ {\em of cellular maps} $C(\sigma) \to X$ {\em extending} $f$
{\em and a number} $\lambda_1 \geq 0$
{\em can be chosen in such a way that for each} $b \in B$ {\em there
  exists} $\tilde{f} \in S$ {\em with}
\[
\gsh_b(\tilde{f}) > \gsh_b(f) - \diam hC(\sigma) - \lambda_1 -
\varepsilon.
\]

\hfill $\Box$
~\\

6.3 {\bf Imposing CAT(0).} We remain in the situation of Section 6.2: We
are given a cellular map $f: C(\mathop\sigma\limits^\bullet) \to X$ which
we try to
extend to $\tilde{f}: C(\sigma) \to X$ keeping $\gsh_a(\tilde{f})$ under
control. But we need control over $\gsh_a(\tilde{f})$ for all points $a
\in M$ and so Corollary 6.3 falls short. In order to improve it we have
to impose the assumption that the metric space $(M,d)$ be CAT(0). The
CAT(0)-condition first appears\footnote{In Remark 7.8 we indicate the other
places in this paper where we use the CAT(0) condition.}
via the following lemma ; see the proof of
[BrHa; III 2.8]:
~\\

{\bf Lemma 6.4.} {\em Let} $(M,d)$ {\em be a} CAT(0) {\em space and
  let} $\varepsilon > 0$. {\em For any} $c \in M$ {\em and} $r \geq 0$,
 {\em any number} $R > r(1+2r/\varepsilon)$ {\em has the property that
   when} $a \in M$ {\em is
  outside} $B_R(c)$ {\em and both} $p$ {\em and} $p'$ {\em are in}
$B_r(c)$ {\em then}
\[
|(d(a,p)-d(a,p')) - (d(b,p) - d(b,p')) | < \varepsilon,
\]
{\em where} $b$ {\em is the point on the geodesic segment from} $c$ {\em
  to} $a$ {\em distant R from c}.

\hfill $\Box$
~\\

We apply Lemma 6.4 in the situation of Corollary 6.3 with the following
careful choice of $B \subseteq M$. Choose a centre $c \in h(\sigma)$ and
a radius $r \geq 0$ such that the ball $B_r(c)$ contains both
$hC(\sigma)$ and $hf(C(\mathop\sigma\limits^\bullet))$. Let $R$ be the radius
given by Lemma 6.4 and put $B:= B_R(c)$.

To improve on Corollary 6.3 we consider a point $a \in M$ outside
$B$. Let $b \in M$ denote the point on the geodesic segment from $c$ to
$a$ with $d(c,b) = R$. Then $b \in B$ so that Corollary 6.3 applies. It
guarantees an extension $\tilde{f} \in S$ with the property that for
each $x \in C(\sigma)$,
\[
\begin{array}{lcl}
d(h\tilde{f}(x),b) & < & d(h(x),b) - \gsh_b(f) + \diam\, hC(\sigma) +
\lambda_1 + \varepsilon\\
& \leq & d(c,b)-\gsh_b(f) + 2\diam\, hC(\sigma) + \lambda_1 + \varepsilon.
\end{array}
\]
Hence
\[ 
\begin{array}{lcl}
d(h\tilde{f}(x),a) & \leq & d(h\tilde{f}(x),b) + d(b,a)\\
& < & d(c,a) - \gsh_b(f) + 2\diam\, hC(\sigma) + \lambda_1 + \varepsilon.
\end{array}
\]
On the other hand,
\[ 
d(h(x),a) \geq d(c,a) - \diam\, hC(\sigma), 
\]
hence
\[d(h(x),a) - d(h\tilde{f}(x),a) > \gsh_b(f) - 3\diam\, hC(\sigma) -
\lambda_1 - \varepsilon.
\]
As the right hand side of this inequality is independent of $x \in
C(\sigma)$ it follows that
\[
\gsh_a(\tilde{f}) \geq \gsh_b(f) - 3\diam\, hC(\sigma) - \lambda_1 -
\varepsilon.
\]
It remains to relate $\gsh_b(f)$ with $\gsh_a(f)$ via Lemma 6.4. For
every $y \in C(\mathop\sigma\limits^\bullet)$ it yields, by putting $p =
h(y)$ and $p' = hf(y)$,
\[
| \sh_{f,a}(y) - \sh_{f,b}(y) |\; < \varepsilon.
\]
It follows that $\gsh_b(f) \geq \gsh_a(f) - \varepsilon$, hence
\[
\gsh_a(\tilde{f}) \geq \gsh_a(f) - 3\diam\, hC(\sigma) - \lambda_1 -
2\varepsilon.
\]
Thus we have proved
~\\

{\bf Proposition 6.5.} {\em Assume that} $X$ {\em is uniformly}
$CC^{n-1}$ {\em over a} CAT(0) {\em space} 
$M$. {\em Then, given any} $\varepsilon > 0$, {\em a
  finite set} $S$ {\em of cellular maps} $C(\sigma) \to X$ {\em
  extending} $f$ {\em and a number} $\lambda_1 \geq 0$ {\em can be chosen in 
such a way that for each} $a \in
M$ {\em there exists} $\tilde{f} \in S$ {\em with}
\[
\gsh_a(\tilde{f}) \geq \gsh_a(f) - 3\diam\, hC(\sigma) - \lambda_1 -
2\varepsilon.
\]

\hfill $\Box$
~\\

6.4 {\bf The main technical theorem.} We are now in a position to prove
the following crucial consequence of the $CC^{n-1}$-condition over a
CAT(0)-space $M$.
~\\

{\bf Proposition 6.6.} {\em Let} $X$ {\em be a cocompact}
$G$-CW-{\em complex with finite stabilizers, let $h : X \to M$ be a control
function and let $n\geq 1$}.  
{\em Assume} $(M,d)$ {\em is} CAT(0) {\em and} $X$ {\em is
  uniformly} $CC^{n-1}$ {\em over} $M$. {\em Then there is a constant}
$\eta \geq 0$ {\em with the property that every homotopically 
closed locally finite}
$G$-{\em sheaf} ${\cal F}: X^{n-1} \leadsto X^{n-1}$ {\em can be
  embedded in a homotopically closed locally finite} $G$-{\em sheaf} $\tilde{{\cal F}}: X^n
\leadsto X^n$ {\em with} $d_a(\tilde{{\cal F}}|\sigma) \leq \eta$ {\em for
  every} $a \in M$ {\em and all} $n$-{\em cells} $\sigma$ {\em of} $X$.
~\\

{\bf Proof.} Fix an $n$-cell $\sigma$ of $X$. As ${\cal F}$ is locally
finite there are only finitely many members $f \in {\cal F}$ with $D(f)
= C(\mathop\sigma\limits^\bullet)$. Hence Proposition 6.5 yields a finite
set $S(\sigma)$ of cellular maps $C(\sigma) \to X^n$ such that
for each $f \in {\cal F}$ with $D(f) = C(\mathop\sigma\limits^\bullet)$
there is some $\tilde{f} \in S(\sigma)$ extending $f$ with $\gsh_a(f) -
\gsh_a(\tilde{f}) \leq \Lambda$, where $\Lambda$ can be chosen
independent of $\sigma$. $S(\sigma)$ depends on $\sigma$; but for
the $G$-translates $g\sigma$ of $\sigma$, we do not have to
choose $S(g\sigma)$ anew but can put $S(g\sigma): = gS(\sigma)$. This
is because $\gsh_{ga}(gf) = \gsh_a(f)$; see \S 5.4. Let $\tilde{{\cal F}}$
denote the $G$-sheaf generated by ${\cal F}$ together with the maps
$\tilde{f} \in S(\tau)$, as $\tau$ runs through a (finite) set of
$n$-cells representing all $G$-orbits. By Proposition 4.4 
$\tilde{{\cal F}}$ fulfills the
conclusion of the Proposition.

\hfill $\Box$
~\\

{\bf Corollary 6.7.} {\em Let} $h: X \to M$ {\em as in Proposition
  6.6. If} $M$ {\em is} CAT(0) {\em and} $X$ {\em is uniformly}
$CC^{n-1}$ {\em over} $M$ {\em then there is a homotopically closed
 locally finite}
$G$-{\em sheaf} ${\cal F}: X^n \leadsto X^n$, {\em a finite radius}
$R \geq 0$, {\em and some} $\varepsilon >0$ {\em such that}
\[
\mu_a({\cal F} | \sigma) - \delta_a({\cal F} | \sigma) \geq
\varepsilon 
\]
{\em for all} $a \in M$ {\em and all cells} $\sigma$ {\em of} $X^n$
  {\em with} $h(\sigma) \cap B_R(a) = \emptyset$.
~\\

{\bf Proof}. Uniformly $CC^{-1}$ means that the
$G$-action on $M$ is cocompact. Hence Theorem 6.1 applies and yields
a radius $R = R(\alpha)$ and a locally finite $G$-sheaf ${\cal F}^0: X^0
\leadsto X^0$ satisfying $\mu_a({\cal F}^0 | v) \geq \alpha$ for each
  vertex $v \in X^0$ with $d(h(v),a) \geq R$, where $\alpha$ is an
  arbitrary positive number yet to be chosen. Now we apply Proposition
  6.6 in each dimension to see that ${\cal F}^0$ can be embedded in a
homotopically closed locally finite $G$-sheaf ${\cal F}$ with a bound on
    $\delta_a({\cal F} | \sigma)$ independent of ${\cal F}^0, \sigma$,
    and $a \in M$. Choosing $\alpha$ greater than or equal to this bound
    yields the corollary.

\hfill $\Box$
~\\

{\bf Theorem 6.8} {\em Let} $h: X \to M$ {\em be a contractible
  G-}CW-{\em complex over M with finite stabilizers and cocompact
  n-skeleton. Under the assumption that the control space M is} CAT(0)
  {\em the following conditions are equivalent}
\begin{enumerate}
\item[(i)] {\em X is uniformly} $CC^{n-1}$ {\em over M.}
\item[(ii)] {\em For every} $a \in M$ {\em there is a G-finitary
    contraction} $\phi_a: X^n \to X^n$ {\em with event radius
    independent of a.}
\item[(iii)] {\em G acts cocompactly on M and there is a G-finitary
    contraction} $\phi: X^n \to X^n$ {\em towards some} $a \in M$.
\item[(iv)] {\em There is a locally finite homotopically closed 
G-sheaf} ${\cal F}: X^n
  \leadsto X^n$, {\em a radius} $R \geq 0$, {\em and a positive number}
  $\varepsilon$ {\em such that}
\[
\mu_a({\cal F}|\sigma) - \delta_a({\cal F}|\sigma) \geq \varepsilon
\]
{\em for all} $a \in M$ {\em and all cells} $\sigma$ {\em of} $X^n$
{\em with} $h(\sigma) \cap B_R(a) = \emptyset$.
\end{enumerate}
~\\

{\em Proof.} As ``uniformly $CC^{-1}$'' implies ``$G$ acts cocompactly on
$M$'' the implication (i) $\Rightarrow$ (iv) is covered by Corollary
6.7. The implication (iv) $\Rightarrow$ (ii) is covered by Proposition
5.7. (ii) $\Rightarrow$ (iii) is easy since the assumption (ii) for $n = 0$
implies, by Proposition 5.2, that $X$ is uniformly $CC^{-1}$ and hence
$M$ is cocompact. Theorem 5.6, when stated for a cocompact $G$-space
$M$, covers the remaining implication (iii) $\Rightarrow$ (i).

\hfill $\Box$
~\\

{\bf Addendum 6.9.}  {\em To the above conditions can be added}:
\begin{enumerate}
\item[(v)] $X$ {\em is uniformly $CC^{n-1}$ over $M$ with constant lag}.
\end{enumerate}

\hfill $\Box$

\newpage

\section{Controlled Connectivity as an Open Condition}

7.1 {\bf The topology on the set of all $\bf G$-actions.} Let $(M,d)$ be a
proper metric space and let $G$ be group of type $F_n$. We consider left
actions of $G$ on $M$ where, in contrast to the previous sections, the
action is allowed to vary. We give $G$ the discrete topology, and we give
$\Isom(M)$ and $\Hom(G,\Isom(M))$ the compact-open topology.

We recall some general topology; for details, see [Du], especially
Section XII 1.3 and 5.2 and IX 9.2. When $A,B$ and $C$ are topological
spaces and $C^B$ denotes the space of all continuous functions $B \to C$
with the compact-open topology, a function $\alpha: A \times B \to C$ is
continuous if and only if its adjoint $\hat{\alpha}: A \to C^B$ is
continuous, provided $B$ is locally compact. (Note that our $M$ is
locally compact.) The space $C^B$ is metrizable provided $B$ and $C$ are
$2^{nd}$ countable, $B$ is locally compact and Hausdorff, and $C$ is
regular. (Note that our $M$ is $2^{nd}$ countable, being metrizable and
 sigma compact.) In the
compact-open topology sequential convergence agrees with uniform
convergence on compact sets ([Du], XII 7.2), provided the target space is
metrizable. Summarizing:
~\\

{\bf Proposition 7.1} $\Isom(M)$ {\em and} $\Hom(G,\Isom(M))$ {\em are
  metrizable} $2^{nd}$ {\em countable spaces. In both function spaces
  the compact-open topology is the topology of uniform convergence on
  compact sets. For any space A, a function} $A \to M^M$ {\em is
  continuous if and only if its adjoint function} $A \times M \to M$
{\em is continuous.}

\hfill $\Box$
~\\

7.2 {\bf Continuous choice of control functions.} In Section 2.2 we
chose a contractible free G-CW complex $X$ with $G$-finite $n$-skeleton
and a $G$ equivariant control function $h: X \to M$ into the given
control space $M$. Now we vary the $G$-action $\rho: G \to \Isom(M)$ and
choose a $(G,\rho)$-equivariant control function $h_\rho: X \to M$ for
each $\rho$. In order to emphasize the dependence on $\rho$ we write
$\rho(g)a$ for the effect of the action of $g \in G$ on $a \in M$. The
assignment $\rho \mapsto h_\rho$ thus defines a function
$\Hom(G,\Isom(M)) \to M^{X^n}$.
~\\

{\bf Proposition 7.2} {\em If (M,d) is a proper unique-geodesic metric
  space this function can be chosen to be
  continuous.}
~\\

{\em Proof.} Pick a base point $b \in M$. Pick a representative cell
$\sigma$ for each $G$-orbit of cells of $X$, and write $x_\sigma \in X$
for the ``barycenter'' of $\sigma$ (i.e., the image of the origin under
a characteristic map $\BB^k \dpo \sigma$). We will construct $h_\rho:
X \to M$ in such a way that $h_\rho(x_\sigma) = b$ for each
representative cell $\sigma$. By $G$-equivariance this defines $h_\rho$
on the zero skeleton $X^0$. Assuming, then, that $h_\rho: X^k \to M$ is
already defined on the $k$-skeleton we extend it to a representative
$(k+1)$-cell $\sigma$ by regarding $\sigma$ as the mapping cone of its
attaching map $S^k \to X^k$. Extend $h_\rho$ to $\sigma$ by mapping each
cone line linearly to the unique geodesic joining the $h_\rho$-images of
its endpoints.  Uniqueness of geodesic implies that geodesics vary
continuously with their endpoints [BrHa I.3.11], so this extension is
continuous. Extend $h_\rho$ to be $(G,\rho)$-equivariant on
$X^{k+1}$. The map $h_\rho: X \to M$ is thus ``canonically defined'' and
is therefore continuous in $\rho$.

\hfill $\Box$
~\\

For the rest of this section we assume that $h_\rho$ depends
continuously on $\rho$. 
We now write $\sh^\rho_{f,a}: D(f) \to \RR$ in order to emphasize
that it depends on the $G$-action. 
~\\

{\bf Proposition 7.3} {\em Let} ${\cal F}: X \leadsto X$ {\em be a
  locally finite homotopically closed G-sheaf on a G-}CW{\em complex X
 over M. Then the
  maximal guaranteed vertex shifts} $\mu^\rho_a({\cal F}|\sigma)$
{\em and the defect} $d^\rho_a({\cal F}|\sigma)$ {\em are jointly
  continuous in the variables} $(\rho,a)$. {\em Moreover
  we have}
\[
|\mu^\rho_a({\cal F}|\sigma) - \mu^\rho_b({\cal F}|\sigma)\,
|\, \leq 2d(a,b)
\]
{\em and}
\[
|d^\rho_a({\cal F}|\sigma) - d^\rho_b({\cal F}|\sigma)\, |\, \leq 4d(a,b)
\]
{\em for all} $a,b \in M$.

\hfill $\Box$
~\\

{\em Proof}.  This is elementary.  It is helpful to use the even more
precise notation $\sh^h_{f,a}(x)$, $\gsh^h_a(f)$ etc. where $h : X \to M$ is
any map, and to prove joint continuity in $h$ and $a$ using Formula (5.7).

\ \hfill $\Box$
~\\

7.3 {\bf Imposing CAT(0).} Condition (iv) of Theorem 6.8 expresses the
$CC^{n-1}$-property of an action $\rho: G \to \Isom(M)$ in terms of the
existence of a locally finite homotopically 
closed $G$-sheaf ${\cal F}: X^n \leadsto
X^n$ whose functions have to satisfy a
certain inequality. In order to prove the ``openness'' of the
$CC^{n-1}$-property we will show that if we fix the $G$-sheaf ${\cal F}$
then this inequality remains true under small perturbation of
the action $\rho$. Proposition 7.3 is the first step in this but is not
sufficient because the range of the parameter $a \in M$ in the
inequality is not compact. The missing ingredient is a result
expressing that the functions of $a, \mu^\rho_a({\cal F}|\sigma)$ and
$d^\rho_a({\cal F}|\sigma)$, are, to some
extent, ruled by their values on points $a \in M$ within a bounded
distance from $h_\rho(\sigma)$. This requires, once again, the assumption that $M$ be CAT(0).

Before we can state the missing result in Proposition 7.5 we need some
preparation. We consider a closed locally finite $G$-sheaf ${\cal F}: X
\leadsto X$ on a $G$-CW complex $X$ over the CAT(0) space $M$. For each
cell $\sigma$ of $X$ we consider the carrier $C(\sigma)$ and the
restricted sheaf ${\cal F}|C(\sigma)$. With respect to a given
$G$-action $\rho: G \to \Isom(M), {\cal F}|C(\sigma)$ has a well defined
norm over $M$ in the sense of \S 5.1. We put $r = r_\sigma(\rho):=
\|{\cal F}|C(\sigma)\| + \diam\, h_\rho C(\sigma)$. Let
\[
(7.1) \quad L_\sigma: \Hom(G,\Isom(M)) \times \RR_{> 0} \to \RR,
\]
be defined by $L_\sigma(\rho,\varepsilon) = r+2r^2/\varepsilon$.
 Observe that joint continuity (resp. $G$-equivariance) of $h_\rho(x)$
 in $\rho$ and $x$ implies that the norm of ${\cal F}|C(\sigma)$ (which
 depends on $\rho$)  and the diameter of
$h_\rho C(\sigma)$ are continuous in $\rho$  (resp. $G$-invariant), hence
~\\

{\bf Proposition 7.4} $L_\sigma$ {\em is continuous, and} $L_{g\sigma}
= L_\sigma$ {\em for every} $g \in G$ {\em and every cell} $\sigma$ {\em
  of X.}

\hfill $\Box$
~\\

In order to keep control over distances between points of $M$ and
$h_\rho$-images of cells we choose a ``barycenter'' $x_\sigma \in
\sigma$ for every cell $\sigma$ of $X$ in such a way that $x_{g\sigma}
= gx_\sigma$ for each $g \in G$. Let $c_\sigma:= h_\rho(x_\sigma)$. Then
$c_{g\sigma} = \rho(g)c_\sigma$ for each $g \in G$.

~\\

{\bf Proposition 7.5} {\em Assume M is} CAT(0). {\em Fix a cell} $\sigma$ {\em
  of X, a number} $\varepsilon > 0$ {\em and a G-action} $\rho: G \to
\Isom(M)$. {\em Then given any radius} $R > L_\sigma(\rho,\varepsilon)$
{\em and a point} $a \in M$ {\em with} $d(c_\sigma,a) \geq R$ {\em we
  have}
\[
(7.2) \quad |\mu^\rho_a({\cal F}|\sigma) - \mu^\rho_b({\cal
  F}|\sigma)\, | \, < \varepsilon,
\]
{\em in general, and}
\[
(7.3) \quad |d^\rho_a({\cal F}|\sigma) - d^\rho_b({\cal F}|\sigma) \,
|\, < 2\varepsilon
\]
{\em when} $\dim\sigma > 0$, {\em where} $b = b_\sigma(\rho,a)$ {\em
  stands for the point on the geodesic segment from} $c_\sigma$ {\em to
  a at distance R from} $c_\sigma$.
~\\

{\em Proof.} The function $r = r_\sigma(\rho)$ used in
(7.1) was chosen so that the ball $B_r(c_\sigma)$ contains both $h_\rho
C(\sigma)$ and $h_\rho f(C(\sigma))$ for each $f \in {\cal
  F}|C(\sigma)$. In this situation we can apply Lemma 6.4 with $p =
h_\rho(x), p' = h_\rho f(x)$ and find, for each $x \in D(f)$,
\[
(7.4) \quad |\sh^\rho_{f,a}(x) - \sh^\rho_{f,b}(x)\, |\, < \varepsilon.
\]
(7.2) is an immediate consequence of (7.4). As to (7.3) we
first observe that (7.4) implies the corresponding inequality,
\[
|\gsh_a(f) - \gsh_b(f)\, |\, < \varepsilon,
\]
for the guaranteed shifts of $f$. Plugging this twice into the
definition of the defect $d_a({\cal F}|\sigma)$ in \S 5.5 yields
the required inequality (7.3).

\hfill $\Box$
~\\

For later use we also record
~\\

{\bf Proposition 7.6} {\em In the situation of Proposition 7.5 the
  points} $b = b_\sigma(\rho,a)$ {\em satisfy the inequality}
\[
d(b_\sigma(\rho,a), b_\sigma(\rho',a)) \leq d(c_\sigma(\rho),
c_\sigma(\rho'))
\]
{\em for any two actions} $\rho, \rho': G \to \Isom(M)$.
~\\

{\em Proof.} By the definition of CAT(0) it suffices to prove the
inequality in the Euclidean plane. Here it is elementary -- e.g. use the
Cosine Theorem to compute the two distances.

\hfill $\Box$
~\\

7.4 {\bf The Openness Theorem.} We are now in a position to prove our
main result.
~\\

{\bf Theorem 7.7} {\em For actions} $\rho$ {\em of a group G of type}
$F_n$ {\em on a proper} CAT(0) {\em space M, uniformly controlled
 n-connectedness is an
  open condition; i.e., if} $\rho$ {\em is uniformly} $CC^{n-1}$ {\em
  over M there is a neighbourhood N of} $\rho$ {\em in}
$\Hom(G,\Isom(M))$ {\em such that every} $\rho' \in N$ {\em is also
  uniformly} $CC^{n-1}$ {\em over M.}
~\\

{\em Proof.} Let $X$ be the universal cover of a
$K(G,1)$-complex with finite $n$-skeleton. Make a continuous choice of
control functions $h_\rho: X^n \to M$, fixing a base point $*$, and of 
barycenters $x_\sigma$ as in the
proof of Proposition 7.2. Let $T$ be a
set of representatives of the $G$-orbits of the cells of $X$. We have
$h_\rho(x_\sigma) = *$ for all $\sigma \in T$. There is a finite set $H
\subseteq G$ such that $H T$ contains all cells of the carriers
$C(\sigma)$ with $\sigma \in T$. 

Now let $\rho$ be a $G$-action on $M$ which is uniformly
$CC^{n-1}$. Let ${\cal F}: X^n \leadsto X^n, R \geq 0$ and $\varepsilon
> 0$ be as given in Condition (iv) of Theorem 6.8. We fix a
natural number $m$ which we will specify later. As we can replace $R$ by
a larger number, if necessary, we may assume that $R > L_\sigma(\rho,
\frac{\varepsilon}{m})$. 

Now we consider the following subsets $N_i(\sigma) \subseteq
\Hom(G,\Isom(M))$.
\[
\begin{array}{rlcl}
&N_1(\sigma) & := & \{\rho'\Big|L_\sigma(\rho',\frac{\varepsilon}{m}) <
R\}\\
&N_2(\sigma) & := & \{\rho'\Big||\mu^\rho_b({\cal F}|\sigma) -
\mu^{\rho'}_b({\cal F}|\sigma)\Big|<\frac{\varepsilon}{m},\, b \in
B_R(*)\}\\ 
&N_3(\sigma) & := & \{\rho'\Big||d^\rho_b({\cal F}|\sigma) -
d^{\rho'}_b({\cal F}|\sigma)\Big|<\frac{\varepsilon}{m},\, b \in 
\rho(H)B_R(*)\}\\
\mbox{and } &N_4 & := & \{\rho'\Big|d(\rho(g)*, \rho'(g)*) < 
\frac{\varepsilon}{m},\, g \in H\}
\end{array}
\]
These sets are open by Propositions 7.3 and 7.4 since $H$ is finite and
$B_R(*)$ is compact. Since $T$ and $HT$ are finite sets of cells the
intersection
\[
N := \bigcap\limits_{\sigma \in T}(N_1(\sigma) \cap N_2(\sigma)) \cap
(\bigcap\limits_{\tau \in HT}N_3(\tau)) \cap N_4
\]
is also open. Of course, $\rho \in N$.

Let $a \in M$ with $d(\rho(g)*,a) \geq R + \frac{\varepsilon}{m}$ for all $g
\in H$. This implies that for each $\tau = g\sigma \in HT$ not
only the element $b_\tau(\rho,a)$ occurring in Proposition 7.5 but also
$b_\tau(\rho',a)$ for all $\rho' \in N_4$ is well defined. Indeed, this
follows since $c_\tau(\rho) = c_{g\sigma}(\rho) = \rho(g) c_\sigma(\rho)
= \rho(g)*$, and similarly for $\rho'$.

For all $\rho' \in N$ and all $\sigma \in T$ we find
\[
|\mu^\rho_a-\mu^{\rho'}_a| \leq |\mu^\rho_a-\mu^{\rho}_b| +
|\mu^\rho_b-\mu^{\rho'}_b| + |\mu^{\rho'}_b-\mu^{\rho'}_a|
\]
\[
< \frac{\varepsilon}{m} + \frac{\varepsilon}{m} + \frac{\varepsilon}{m}
= \frac{3\varepsilon}{m}.
\]
Here we have omitted the argument ${\cal F}|\sigma$ and we have
used $b$ to denote the point $b = b_\sigma(\rho,a)$; we have used
$c_\sigma(\rho) = * = c_\sigma(\rho')$ so that $b_\sigma(\rho,a) =
b_\sigma(\rho',a)$ and have applied Proposition 7.5 both to $\rho$ and
$\rho'$; and we have used $b \in B_R(*)$ so that the definition of
$N_2(\sigma)$ applies.

The corresponding computation for $d^\rho_a$ instead of $\mu^\rho_a$ is
more subtle since we will need it not only for $\sigma \in T$ but for
all $\tau = g\sigma \in HT$ in order to get control over the total
defect $\delta^\rho_a$. Again omitting the argument ${\cal F}|g\sigma$
we find the inequality
\[
|d^\rho_a-d^{\rho'}_a| \leq |d^\rho_a-d^\rho_b| +
|d^\rho_b-d^{\rho'}_{b'}| + |d^{\rho'}_{b'}-d^{\rho'}_a|,
\]
where $b = b_{g\sigma}(\rho,a)$ and $b' = b_{g,\sigma}(\rho',a)$. By
Proposition 7.5 the first and third terms of the right hand side are
smaller than $\frac{2\varepsilon}{m}$. But the middle term needs more
care. We write
\[
|d^\rho_b-d^{\rho'}_{b'}| \leq |d^\rho_b-d^{\rho'}_b| + |
d^{\rho'}_b-d^{\rho'}_{b'}|
\]
and observe that $b = b_{g\sigma}(\rho,a) \in B_R(c_{g\sigma}(\rho)) =
\rho(g)B_R(*)$. Hence the definition of $N_3(g\sigma)$ applies and shows
that the first terms of the right hand side is less than
$\frac{\varepsilon}{m}$. The second term is at most equal to $4d(b,b')$
by Proposition 7.3 which is less than or equal to $4d(c_{g\sigma}(\rho),
c_{g\sigma}(\rho')) = 4d(\rho(g)*, \rho'(g)*) \leq
  \frac{4\varepsilon}{m}$, by Proposition 7.6.

Thus we have proved that for the defect on a cell $\tau \in HT$ we
have
\[
|d^\rho_a({\cal F}|\tau) - d^{\rho'}_a({\cal F}|\tau) \Big| <
\frac{9\varepsilon}{m}.
\]
This is enough information to obtain a similar formula for the total
defect on the cells $\sigma \in T$. Indeed, we find
\[
|\delta^\rho_a({\cal F}|\sigma) - \delta^{\rho'}_a({\cal F}|\sigma)
\Big| < \frac{9\varepsilon \cdot \dim\sigma}{m} \leq \frac{9\varepsilon
  n}{m}
\]
Putting things together we find that if $\rho' \in N$ then we have on
each cell $\sigma \in T$
\[
|(\mu^\rho_a-\delta^\rho_a) - (\mu^{\rho'}_a-\delta^{\rho'}_a) \Big|
< \frac{3(3n+1)\varepsilon}{m}.
\]
Now choose $m := 6(3n+1)$ so that, since $\rho$ satisfies (iv) of
Theorem 6.8,
\[
|\mu^{\rho'}_a({\cal F}|\sigma) - \delta^{\rho'}_a({\cal F}|\sigma)
\Big| > \frac{\varepsilon}{2}
\]
for all $\sigma \in T$. Our computation applies for all $a \in M$ with
$d(c_\sigma(\rho'),a) = d(*,a) \geq R_0$, where $R_0 =
R+\frac{\varepsilon}{m}+\max d(\rho(F)*,*)$. All cells are translates of
cells in $T$, so that the $G$-compatibility conditions in \S\S 5.4 and
5.5 imply that condition (iv) of Theorem 6.8 holds true for $\rho'$.

\hfill $\Box$
~\\

{\bf Remark 7.8.} The case $n = 0$ -- i.e. the condition that $M$ be
cocompact -- is of course simpler. First, it requires
Theorem 6.1 rather than Theorem 6.8, and Theorem 6.1 holds true
without the assumption that $M$ be CAT(0). Secondly, the defect function
$d^\rho_a$ plays no role in dimension 0, so the more subtle
technicalities of the above proof disappear -- in particular Proposition
7.6, which was another instance where the CAT(0) condition appeared, is
not needed. However, the case $n = 0$ does require the first part of
 Proposition 7.5 which was the third instance where CAT(0) was needed.
 We do not know to what extent an openness theorem for cocompact
 actions on geodesic metric spaces beyond the CAT(0) case holds.

\newpage
\section{Completion of the proofs of Theorems A and A$'$}

8.1 {\bf Controlled acyclicity.}  In the set-up of \S2.3, we say $X$ is {\em
controlled} $(n-1)$-{\em acyclic} $(CA^{n-1})$ over $a$ (with respect to
$h$) if for all $r \geq 0$ and $-1\leq p\leq n-1$ there exists $\lambda \geq
0$ such that every $\ZZ$-cycle in $X_{(a,r)}$ bounds in
$X_{(a,r+\lambda)}$.
~\\

8.2.  {\bf The $F_n$-Criterion.}  In [Br 87$_{\mbox{I}}$, Theorem 2.2] Brown
gives an $FP_n$-Criterion analogous to the $F_n$-Criterion:  i.e. $F_n$ is
replaced by $FP_n$ and $CC^{n-1}$ is replaced\footnote{In Brown's
terminology the corresponding homotopy or homology is ``essentially
trivial'' in dimensions $\leq n-1$.} by $CA^{n-1}$.   From this, Brown
deduces the ``if'' direction of the $F_n$-Criterion by proving it for $n\leq
2$:  this together with the $FP_n$-Criterion gives ``$F_2$ and $FP_n$''
which is well known to be equivalent to $F_n$.  He also proves the ``only
if'' direction of the $F_2$-Criterion. To complete the proof of the ``only
if'' direction of the $F_n$-Criterion via the $FP_n$-Criterion one needs a
``Hurewicz-type'' theorem for which we do not know a reference, but which
will appear in [Ge]: 
~\\

{\bf Theorem 8.1.} {\em Let} $(Y_r)_{r\in \RR}$ {\em be a filtration of the}
CW {\em complex} $Y$ {\em by subcomplexes, and let} $f : Y \to \RR$ {\em be 
the control function} $f(y) = \inf\{r\mid y\in Y_r\}$.  {\em With respect to}
$f$, $Y$ {\em is} $CC^{n-1}$ {\em if and only if} $Y$ {\em is} $CC^1$ 
{\em and} $CA^{n-1}$.  

\hfill $\Box$
~\\

8.3 {\bf Proof of Theorem A.}  It only remains to prove:
~\\

{\bf Lemma 8.2.}  {\em The} $G_a$-{\em complexes} $X^n_{(a,r)}$ {\em are 
cocompact}.
~\\

{\em Proof}.  The cells of $X^n_{(a,r)}$ decompose into finitely many 
classes under
the $G$-action, and those are, of course, $G_a$-invariant. 
We claim that for every cell $\sigma$ of $X^n_{(a,r)}$ the set of all
cells of the form $g\sigma \in X^n_{(a,r)},\, g \in G$, consists of only
finitely many $G_a$-orbits.

For every cell $\sigma$
of $X^n$ we choose  a ``barycenter'' $x_\sigma \in \sigma$ in such a way that
$x_{g\sigma} = gx_\sigma$ for all $g \in G$. Since $G$ acts discretely
on $M$ and $X^n$ has only finitely many $G$-orbits of cells, $B:=
\{h(x_\sigma)|\sigma \subset X^n\}$ is a discrete subset of $M$. Since
$B_r(a)$ is compact the subset $C:= B \cap B_r(a)$ is finite. For each
$c \in C$ the stabilizers $G_a$  and $G_c$ are commensurable and we
choose a (finite) set of coset representatives for $(G_a \cap
G_c)\backslash G_c$. Putting these together, we get a finite
subset $L \subset G$ with $G_c \subset G_aL$, for all $c \in
C$. Moreover, we consider a second finite subset $T \subset G$ by
choosing, for each pair of elements $c_1,c_2 \in C$ which are in the
same $G$-orbit an element $t \in G$ with $g_2 = tg_1$.

We complete the proof of the Claim by showing that if $g \in G$ and
$\sigma$ is a cell of $X^n$ such that both $\sigma$ and $g\sigma$ are in
$X_{(a,r)}$ then $g \in G_aLT$. Indeed, both $c:= h(x_\sigma)$ and $gc =
h(x_{g\sigma})$ are in $C$; hence there is some $t \in T$ with $gc =
tc$. This shows that $gt^{-1} \in G_{tc}$, hence $g \in G_{tc}T$. But
$tc \in C$, so that $G_{tc} \subseteq G_aL$. This proves the Claim.

\hfill $\Box$
~\\

8.4 {\bf Properly discontinuous actions.} Let $(M,d)$ be a metric space
and $Q$ a group which acts on $M$ by isometries. Recall that this action
is properly discontinuous if every point of $M$ has a neighbourhood $U$
such that $\{q \in Q | U \cap qU \not= \emptyset\}$ is finite. To
complete Footnote 8 (and hence the proof that Theorem A implies Theorem
A$'$)we only need:
~\\

{\bf Lemma 8.5} {\em If M is locally compact then the Q-action on M is
  properly discontinuous if and only if its orbits are discrete and its  
point stabilizers are finite.}
~\\

{\em Proof.} For each $a \in M$ and $\varepsilon > 0$ we consider the
subset $L(a,\varepsilon) \subset Q,\, L(a,\varepsilon):= \{q | d(a,qa) <
\varepsilon\}$. Observe that $L(a,\varepsilon)$ contains the stabilizer
$Q_a$. On the one hand, to say that the orbit $Q_a$ is a closed discrete
subspace of $M$ is to say that there is some $\varepsilon >
0$ with $L(a,\varepsilon) = Q_a$. On the other hand, the 
definition of ``properly discontinuous'', when translated into a statement
about the metric $d$ , says that for every $a \in M$ there is some
$\varepsilon > 0$ such that $L(a,\varepsilon)$ is finite. But then one
can take $0 < \nu < \varepsilon$ with $d(a,qa) > \nu$, for all $q \in
L(a,\varepsilon)$ with $qa \not= a$, and one has $Q_a = L(a,\nu)$ finite
for some $\nu > 0$. This proves the Lemma.

\hfill $\Box$

\newpage
\section{The Invariance Theorem}

Here we prove the Invariance Theorem 3.3 which says that in checking the
$CC^{n-1}$ property over $a \in M$ we may use any $a \in M$, any
$n$-dimensional, $(n-1)$-connected, cocompact rigid $G$-CW complex $Y$ 
and any $G$-map $h: Y \to M$, provided only that the stabilizer of each
$p$-cell of $Y$ is of type $F_{n-p}, 0 \leq p \leq n-1$.  Recall that in \S3
we have already proved Invariance among such complexes $Y$ which are free:
we will use that in the present proof.  

Independence of $a$ is clear. Given $Y$, independence of $h$ is proved as
in Proposition 3.1. To establish independence of $Y$, we let $X$ be a
contractible free $G$-CW complex. Then $X \times Y$ is an
$(n-1)$-connected free $G$-CW complex. We form a commutative diagram:
\[
\begin{array}{lcccl}
\tilde{W} & \stackrel{\tilde{\alpha}}{\longrightarrow} & X \times Y &
\stackrel{\tilde{\beta}}{\longrightarrow} & Y 
\stackrel{h}{\longrightarrow} M\\
\downarrow & & \downarrow & & \downarrow\\
W & \stackrel{\alpha}{\longrightarrow} & G\backslash(X\times Y) &
\stackrel{\beta}{\longrightarrow} & G\backslash Y
\end{array}
\]
The space $G\backslash (X\x Y)$ is the $n$-skeleton of a $K(G,1)$-complex.
The $G$-map $\tilde\beta$ is projection on the $Y$-factor.  There is a
cellular map $\beta$ making the square involving $\beta$ and $\tilde\beta$
(in the diagram) commute.

Consider a cell $\sigma$ of $G\backslash Y$.  By rigidity it is the
homeomorphic image of a cell of $Y$, and if $z\in \mathop\sigma\limits^\circ$ 
then $Z(\sigma) := \beta^{-1}(z)$ is the $n$-skeleton of a
$K(G_\sigma,1$)-complex (though not a subcomplex of $G\backslash (X\x Y)$
when dim $\sigma > 0$).  Indeed, $G\backslash (X\x Y)$ can be
decomposed\footnote{Compare the ``graphs of complexes'' in [SW; p. 165 et
seq.] which the present argument generalizes, and see [Ge, Ch 2] or [Sta,
\S3] for details omitted here.}  as a ``CW complex of CW complexes'',
parametrized by the cells $\sigma$ of $G\backslash Y$, in
which each complex $C(\sigma) := Z(\sigma) \x B^{{\rm{dim}}
(\sigma)}$ is
glued to $\beta^{-1}(G\backslash Y^{{\rm{dim}}(\sigma)-1})$ 
along $Z(\sigma)
\x S^{{\rm{dim}}(\sigma)-1}$.  It follows that, for each $y\in Y$,
$\tilde\beta^{-1}(y)$ an $n$-dimensional $(n-1)$-connected $G_y$-complex.  

Replacing the complex $Z(\sigma)$ by a new complex $Z'(\sigma)$ of the same
homotopy type but with finite $(n-\mbox{dim}(\sigma))$-skeleton, one gets
$W$, again the $n$-skeleton of a $K(G,1)$-complex, by gluing the complexes
$C'(\sigma) := Z'(\sigma) \x B^{{\rm{dim}}(\sigma)}$, as before, by
induction on dim $\sigma$.  The map $\alpha$ naturally arising in this
construction takes each $C'(\sigma)$ to $C(\sigma)$.  Since $G\backslash Y$
is a finite complex, one sees by counting cells that $W$ is a finite
complex.  The commutative diagram is completed with the universal cover
$\tilde W$ of $W$ and the lift $\tilde\alpha$ of $\alpha$.  The map
$\tilde\alpha$ is a $G$-map, and for each $y\in Y$
$\tilde\alpha^{-1}\tilde\beta^{-1}(y)$ is $n$-dimensional and
$(n-1)$-connected.  It follows by the Vietoris-Smale Theorem [Sm 57] that
for every subcomplex $K$ of $Y$, $\tilde\beta\circ \tilde\alpha \mid :
\tilde\alpha^{-1}\tilde\beta^{-1}(K) \to K$ induces an isomorphism on
homotopy groups in dimensions $\leq n-1$.  Thus $Y$ is $CC^{n-1}$ over $a$
with respect to the control function $h$ if and only $\tilde W$ is $CC^{n-1}$
over $a$ with respect to $h\circ \tilde\beta \circ \tilde \alpha$.  But
$\tilde W$ is a free $G$-complex. Since invariance has already been proved
for free $G$-complexes, we are done.  

\hfill $\Box$

\newpage

\fontsize{10}{10pt} \selectfont  
\noindent {\large{\bf References}}
~\\

\bigspitem{BeBr 97} M. Bestvina and N. Brady, Morse theory and
finiteness properties of groups, {\em Invent. Math.} {\bf 129} (1997),
445-470.

\bigspitem{BG$_{\mbox{II}}$} R. Bieri and R. Geoghegan, Connectivity
properties of group actions on non-positively curves spaces II:  the
geometric invariants, preprint.

\bigspitem{BGr 84} R. Bieri and J.R.J. Groves, The geometry of the set of
characters induced by valuations, {\em J. reine und angew. Math.} {\bf 347}
(1984), 168-195.

\bigspitem{BNS 87} R. Bieri, W. Neumann and R. Strebel, A geometric
invariant of discrete groups, {\em Invent. Math.} {\bf 90} (1987), 451-477.

\bigspitem{BRe 88} R. Bieri and B. Renz, Valuations on free resolutions and
higher geometric invariants of groups, {\em Comment Math. Helvetici} {\bf
63} (1988), 464-497.

\bigspitem{BS 80} R. Bieri and R. Strebel, Valuations and finitely presented
metabelian groups, {\em Proc. London Math. Soc.} (3) {\bf 41} (1980),
439-464.

\bigspitem{BS} R. Bieri and R. Strebel, {\em Geometric invariants for
discrete groups}, (monograph in preparation).

\bigspitem{BrHa} M. Bridson and A. Haefliger, {\em Metric spaces of non-positive
curvature}, (monograph in preparation).

\bigspitem{Br 87$_{\mbox{I}}$} K.S. Brown, Finiteness properties of groups,
{\em J. Pure and Applied Algebra} {\bf 44} (1987), 45-75.

\bigspitem{Bu} K.-U. Bux, {\em Endlichkeitseigenschaften aufl{\"o}sbarer
$S$-arithmetischer Gruppen {\"u}ber Funktionenk{\"o}rpern}, Dissertation,
Frankfurt, 1998.

\bigspitem{Du} J. Dugundji, {\em Topology}, Allyn and Bacon, Boston, 1966.

\bigspitem{Fa 99} F.T. Farrell, Fibered representations, an open condition,
{\em Topology and its Applications}, to appear.  

\bigspitem{FePe 95} S. Ferry and E.K. Pedersen, Epsilon surgery theory,
{\em Proceedings 1993 Oberwolfach conference on Novikov conjectures, 
rigidity and
index theorems}, vol. 2, (A. Ranicki, ed.) 
London Math. Soc. Lecture Notes, vol. 227, Cambridge
University Press (1995), 167-226.

\bigspitem{FrLe 85} D. Fried and R. Lee, Realizing group automorphisms, {\em
Contemp. Math.} {\bf 36} (1985), 427-433.

\bigspitem{Ge} R. Geoghegan, {\em Topological methods in group theory}, 
(monograph in preparation).

\bigspitem{Geh} R. Gehrke, {\em Die h{\"o}hern geometrischen Invarianten
f{\"u}r Gruppen mit Kommutatorrelationen}, Dissertation, Frankfurt, 1992.

\bigspitem{Hu} S.-T. Hu, {\em Theory of retracts}, Wayne State University
Press, Detroit, 1965.

\bigspitem{Ko 96} D. Kochloukova, The FP$_m$-conjecture for a class of
metabelian groups, {\em J. of Algebra} {\bf 184} (1996), 1175-1204.

\bigspitem{Me 94} H. Meinert, The geometric invariants of direct products of
virtually free groups, {\em Comment. Math. Helvetici} {\bf 69} (1994),
39-48.

\bigspitem{Me 95} H. Meinert, The Bieri-Neumann-Strebel invariant for graph
products of groups, {\em J. Pure Appl. Algebra} {\bf 103} (1995), 205-210.

\bigspitem{Me 96} H. Meinert, The homological invariants of metabelian
groups of finite Pr{\"u}fer rank:  a proof of the $\Sigma^m$-conjecture,
{\em Proc. London Math. Soc.} (3) {\bf 72} (1996), 385-424.

\bigspitem{Me 97} H. Meinert, Actions on 2-complexes and the homotopical
invariant $\Sigma^2$ of a group, {\em J. Pure Appl. Algebra}, {\bf 119}
(1997), 297-317.

\bigspitem{MMV 98} J. Meier, H. Meinert, and L. VanWyk, Higher generation 
subgroup sets and the $\Sigma$-invariants of graph groups, {\em Comment.
Math. Helv.} {\bf 73} (1998), 22-44.

\bigspitem{Ne 79}  W.D. Neumann, Normal subgroups with infinite cyclic
quotient, {\em Math. Sci.} {\bf 4} (1979), 143-148.

\bigspitem{Ra} M.S. Ragunathan, {\em Discrete subgroups of Lie groups}, 
Ergebnisse
der Math., vol. 68, Spring, Berlin, 1972.

\bigspitem{Re 88} B. Renz, {\em Geometrische Invarianten und
Endlichkeitseigenschaften von Gruppen}, Dissertation, Frankfurt (1988).

\bigspitem{Re 89} B. Renz, Geometric invariants and HNN-extensions, {\em Group
Theory (Singapore 1987)}, de Gruyter Verlag, Berlin 1989, 465-484.

\bigspitem{SW 79} G.P. Scott and C.T.C. Wall, Topological methods in group
theory, {\em Homological group theory} (C.T.C. Wall, ed.) London Math. Soc.
Lecture Notes 36, Cambridge University Press, Cambridge, 1979, 137-203.

\bigspitem{Se} J.-P. Serre, {\em Trees}, Springer, Berlin, 1980.

\bigspitem{Sm 57} S. Smale, A Vietoris mapping theorem for homotopy, {\em
Proc. Amer. Math. Soc.} {\bf 8} (1957), 604-610.

\bigspitem{Sta 98} C. Stark, Approximate finiteness properties of infinite
groups, {\em Topology} {\bf 37} (1998) 1073-1088.
 
\bigspitem{St 80} U. Stuhler, Homological properties of certain arithmetic
groups in the function field case, {\em Invent. Math.} {\bf 57} (1980),
263-281.

\bigspitem{We 64} A. Weil, Remarks on the cohomology of groups, {\em Ann.
Math.} {\bf 80} (1964), 149-157.
~\\

\noindent \textsc{Robert Bieri, Fachbereich Mathematik, Robert-Mayer-Strasse 
6-10, 60325 Frankfurt, Germany}

\noindent {\em e-mail}:  bieri@math.uni-frankfurt.de
~\\

\noindent \textsc{Ross Geoghegan, Department of Mathematics, 
State University of New York, Binghamton, NY  13902-6000, USA}

\noindent {\em e-mail}:  ross@math.binghamton.edu

\end{document}